\documentclass{article} 

\usepackage[english]{babel} 
\usepackage{amssymb}
\usepackage{amsmath}
\usepackage{txfonts}
\usepackage{mathdots}
\usepackage[classicReIm]{kpfonts}
\usepackage{graphicx}
\usepackage{subcaption}

\usepackage[dvipsnames]{xcolor}

\begin{document}


\noindent 

\noindent 

\noindent 

\noindent 

\begin{center}
\textbf{{\Large An exact algorithm for linear optimization problem subject to max-product fuzzy relational inequalities with fuzzy constraints}}
\end{center}

\noindent 
\begin{center}
\noindent \textbf{Amin Ghodousian${}$${}^{*}$\footnote{*\ Corresponding\ author$\\ 
$ Email\ addresses:\ \ $  $a.ghodousian@ut.ac.ir\ (A.\ Ghodousian).}, Romina Omidi${}^{2}$}
\end{center}
\noindent ${}^{1\ }${{\small Faculty of Engineering Science, College of Engineering, University of Tehran, P.O.Box 11365-4563, Tehran, Iran.}}

\noindent ${}^{2}$${}^{\ }${{\small Department of Electrical and Computer Engineering, College of Engineering, University of Tehran, Tehran, Iran.}}

\vskip 0.2in
\noindent \textbf{}

\noindent \textbf{Abstract}
\vskip 0.2in
\noindent Fuzzy relational inequalities with fuzzy constraints (FRI-FC) are the generalized form of fuzzy relational inequalities (FRI) in which fuzzy inequality replaces ordinary inequality in the constraints. Fuzzy constraints enable us to attain optimal points (called super-optima) that are better solutions than those resulted from the resolution of the similar problems with ordinary inequality constraints. This paper considers the linear objective function optimization with respect to max-product FRI-FC problems. It is proved that there is a set of optimization problems equivalent to the primal problem. Based on the algebraic structure of the primal problem and its equivalent forms, some simplification operations are presented to convert the main problem into a more simplified one. Finally, by some appropriate mathematical manipulations, the main problem is transformed into an optimization model whose constraints are linear. The proposed linearization method not only provides a super-optimum (that is better solution than ordinary feasible optimal solutions) but also finds the best super-optimum for the main problem. The current approach is compared with our previous work and some well-known heuristic algorithms by applying them to random test problems in different sizes.

\vskip 0.2in

\noindent \textbf{}

\noindent \textbf{1. Introduction}

\vskip 0.2in
\noindent The theory of fuzzy relational equations (FRE) as a generalized version of Boolean relation equations was firstly proposed by Sanchez and was applied to problems related to the medical diagnosis [41]. Pedrycz categorized and extended two ways of the generalizations of FRE in terms of sets under discussion and various operations which are taken account [36]. Since then, FRE was applied in many other fields such as fuzzy control, prediction of fuzzy systems, fuzzy decision making, fuzzy pattern recognition, image compression and reconstruction, fuzzy clustering and so on. Generally, when rules of inference are applied and their corresponding consequences are known, the problem of determining antecedents is simplified and mathematically reduced to solving an FRE [34]. Nowadays, it is well known that many of the issues associated to the body knowledge can be treated as FRE problems [35].

\noindent The solvability identification and finding set of solutions are the primary, and the most fundamental, matter concerning the FRE problems. Di Nola et al. proved that the solution set of FRE (if it is nonempty), defined by continuous max-t-norm composition is often a non-convex set. This non-convex set is completely determined by one maximum solution and a finite number of minimal solutions [5]. Such non-convexity property is one of two bottlenecks making a major contribution towards an increase in complexity of FRE-related problems, particularly, in the optimization problems subjected to a system of fuzzy relations. Another bottleneck point is concerned with detecting the minimal solutions for FREs. Chen and Wang [2] presented an algorithm for obtaining the logical representation of all minimal solutions and deduced that a polynomial-time algorithm with the ability to find all minimal solutions of FRE (with max-min composition) may not exist. Also, Markovskii showed that solving max-product FRE is closely related to the covering problem which is a type of NP-hard problem [33].  In fact, the same result holds true for a more general t-norms instead of the minimum and product operators [3,29,30]. Over the past decades, the solvability of FRE which is defined using different max-t compositions have been investigated by many researchers [15,16,19,39,43,46,47,50,52]. Moreover, some other researchers have worked on introducing novel concept and at times improving some of the existing theoretical aspects and applications of fuzzy relational inequalities (FRI) [13,17,18,20,21,27,53]. Li and Yang [27] studied FRI with addition-min composition and presented an algorithm to search for minimal solutions. They applied FRI to data transmission mechanism in a BitTorrent-like Peer-to-Peer file sharing systems. In [13], the authors focused on the algebraic structure of two fuzzy relational inequalities \(A \varphi x \leq b^{1}\) and \(D \varphi x \geq b^{2}\). Their research focuses on the study of a mixed fuzzy system formed by two of the earlier FRIs, where \(\varphi\) is an operator with (closed) convex solutions. Generally, if \(\varphi\) is an operator with closed convex solutions, the set of solutions for \(D \varphi x \geq b^{2}\) is determined by a finite number of maximal solutions as well as the same number for minimal ones. In particular, if \(\varphi\) is a continuous non-decreasing function (specially, a continuous t-norm), all maximal solutions overlap each other [13]. Guo et al. [21] investigated a type of FRI problems and the relationship between minimal solutions and FRI paths. They also introduced some rules for reducing the problems and presented an algorithm for solving optimization problems using FRI constraints.

\noindent The optimization problem subject to FRE and FRI is one of the most interesting and on-going research topics amongst similar problems [1,8,12,13,15-23,26,31,40,44,48,53]. Many methods were designed based on the translation of the main problem into an integer linear programming problem which is then solved using well-developed techniques. On the contrary, other algorithms benefit the resolution of the feasible region, some necessary and sufficient conditions for the optimality and simplification processes. The most methods of this category are based on analytical results provided mainly by Sanchez [42] and Pedrycz [37]. Fang and Li converted a linear optimization problem subjected to FRE constraints with max-min operation into an integer programming problem and solved it by a branch-and-bound method using jump-tracking technique [9].  Wu et al. worked on improvement of the method employed by Fang and Li; this was done by decreasing the search domain and presented a simplification process by three rules which resulted from a necessary condition [49]. Chang and Shieh presented new theoretical results concerning the linear optimization problem constrained by fuzzy max–min relation equations [1]. They improved an upper bound on the optimal objective value, some rules for simplifying the problem and proposed a rule for reducing the solution tree. In [25], an application of optimizing the linear objective with max-min composition was employed for the streaming media provider seeking a minimum cost while fulfilling the requirements assumed by a three-tier framework. Linear optimization problem was further investigated by numerous scholars focusing on max-product operation [23,32]. Loetamonphong and Fang defined two sub-problems by separating negative from non-negative coefficients in the objective function, and then obtained an optimal solution by combining the optimal solutions of the two sub-problems [32]. The maximum solution of FRE is the optimum for the sub-problem having negative coefficients. Another sub-problem was converted into a binary programming problem and was solved using a branch-and-bound method. Also, in [23] some necessary conditions to test the feasibility and simplification techniques were presented in order to solve FRE with max-product composition. Moreover, generalizations of the linear optimization problem with respect to FRE have been studied; this was done through replacement of max-min and max-product compositions with different fuzzy compositions such as max-average composition [48] or max-t-norm composition [15,16,19,22,26,44]. For example, Li and Fang solved the linear optimization problem subjected to a system of sup-t equations by reducing it to a 0-1 integer optimization problem [26]. In [22], a method was presented for solving linear optimization problems with the max-Archimedean t-norm fuzzy relation equation constraint. In [44], the authors solved the same problem whit continuous Archimedean t-norm, and to obtain some optimal variables, they used the covering problem rather than the branch-and-bound methods.     

\noindent Recently, many interesting forms of generalizations of the linear programming applied to the system of fuzzy relations have been introduced, and developed based on composite operations used in FRE, fuzzy relations used in the definition of the constraints, some developments on the objective function of the problems and other ideas [4,6,10,28,31,51].  For example, Wu et al. represented an efficient method to optimize a linear fractional programming problem under FRE with max-Archimedean t-norm composition [51]. Dempe and Ruziyeva generalized the fuzzy linear optimization problem by considering fuzzy coefficients [4]. In addition, Dubey et al. studied linear programming problems involving interval uncertainty modeled using intuitionistic fuzzy set [6]. The linear optimization of bipolar FRE was also the focus of study carried out by some researchers where FRE was defined with max-min composition [10] and max-Lukasiewicz composition [28,31]. For example, in [28], the authors introduced a linear optimization problem subjected to a system of bipolar FRE defined as \(X(A^{+}, A^{-}, b)=\{x \in [0,1]^{m}: x \circ A^{+} \vee  \tilde{x} \circ A^{-} = b\}\) where \(\tilde{x}_{i}=1-x_{i}\) for each component of \(\tilde{x}=(\tilde{x}_{i})_{1 \times m}\) and the notations “\(\vee\)” and “\(\circ\)” denote max operation and the max-Lukasiewicz composition, respectively. They translated the original problem into a 0-1 integer linear problem which is then solved using well-developed techniques. In a separate, the foregoing bipolar linear optimization problem was solved by an analytical method based on the resolution and some structural properties of the feasible region (using a necessary condition for characterizing an optimal solution and a simplification process for reducing the problem) [31].

\noindent The optimization problem subjected to various versions of FRI is widely available in the literature as well [12,13,17,18,20,21,53,54]. In [13], the authors studied the linear optimization with constraints formed by \(X(A, D, b^{1}, b^{2}) = \{x \in [0,1]^{n}: A \varphi x \leq b^{1}, D \varphi x \geq b^{2}\}\) where \(\varphi\) represents an operator with convex solutions (e.g., non-decreasing or non-increasing operator). They showed that the feasible region can be expressed as the union of a finite number of convex sets. In particular, if \(\varphi\) is an operator with closed convex solutions such as continuous non-decreasing (non-increasing) operator, the preceding convex sets become closed as well. Therefore, since each t-norm is a non-decreasing function (resulted directly from the property of the monotonicity in t-norms), continuous t-norms introduce important especial examples of operators with closed convex sets. For this reason, the authors proved that the feasible solutions set defined by a continuous t-norm can be formed as the union of a finite number of compact convex sets. Additionally, because of the identity law of t-norms, it was concluded that the feasible region actually consists of points being between one maximum solution and a finite number of minimal solutions. Yang studied the optimal solution of minimizing a linear objective function subject to a FRI where the constraints defined as \(a_{i1}  \wedge x_{1} + a_{i2}  \wedge x_{2} + ... + a_{in}  \wedge x_{n} \geq b_{i}\) for \(i=1,...,m\) and \(a \wedge b = min\{a,b\}\) [54]. He presented an algorithm based on some properties of the minimal solutions of the FRI.  Also, in [53], the authors introduced the latticized linear programming problem subject to max-product fuzzy relation inequalities with application in the optimization management model for wireless communication emission base stations. The latticized linear programming problem was defined by minimizing the objective function \(z(x)=x_{1} \vee x_{2} \vee ... \vee x_{n}\) subject to the feasible region \(X(A, b) = \{x \in [0,1]^{n}:A \circ x \geq b\}\) where “\(\circ\)” denotes fuzzy max-product composition. They also presented an algorithm based on the resolution of the feasible region.

\noindent The FRI-FC problem was introduced in [12] as the following mathematical model in which \(\varphi\) is the minimum t-norm:

\begin{equation} \label{GrindEQ__1_} 
\begin{array}{l} {\min \, \, \, \, \, c^{T}x} \\ {\, \, \, \, \, \, \, \, \, \, \, \, \,A \varphi x \circ b} \\ {\, \, \, \, \, \, \, \, \, \, \, x \in [0,1]^{n}} \end{array} 
\end{equation}

\noindent where \(A=(a_{ij})_{m \times n}\) is a fuzzy matrix and \(b=(b_{i})_{m \times 1}\) is a fuzzy vector such that \(0 \leq a_{ij} \leq 1\) and \(0 \leq b_{i} \leq 1\) for each \(i \in I = \{1,2,...,m\}\) and each \(j \in J = \{1,2,...,n\}\), the constant vector \(c=(c_{j})_{n \times 1}\) and the unknown vector \(x=(x_{j})_{n \times 1}\) are in \(\mathbb{R}^{n}\), \(A \varphi x \circ b\) denotes a fuzzy max-\(\varphi\) composition and “\(\circ\)” denotes relaxed or fuzzy version of the ordinary inequality “\(\leq\)”. If \(a_{i}\) denotes \(i\)'th row of matrix \(A\), then \(i\)'th constraint of problem (1) can be expressed as \(a_{i} \varphi x \circ b_{i}\) which means \(\underset{j=1}{\overset{n}{max}}\{\varphi(a_{ij}, x_{j})\} \circ b_{i}\), \(\forall i \in I\). So, problem (1) can be also interpreted as a generalization of the following problem with the ordinary inequality [13]:

\begin{equation} \label{GrindEQ__1_} 
\begin{array}{l} {\min \, \, \, \, \, c^{T}x} \\ {\, \, \, \, \, \, \, \, \, \, \, \, \,A \varphi x \leq b} \\ {\, \, \, \, \, \, \, \, \, \, \, x\in [0,1]^{n} } \end{array} 
\end{equation}

\noindent In [12], the authors presented an analytical constructive algorithm to find a super-optimum for problem (1), and used this fuzzy system to convincingly optimize the quality of education in school education (with minimum cost) while were to be selected by parents. In [14], the authors presented a modified PSO algorithm for solving the most general version of problem (1) where \(\varphi\) is an arbitrary continuous t-norm. They showed that the modified PSO can produce high quality solutions with low average error (good accuracy) and low standard deviation (good stability). Also, it was shown that the modified PSO generates better solutions compared with the solutions obtained by the max-min algorithm [12] and the original PSO [24] when \(\varphi\) is considered as the minimum t-norm and an arbitrary t-norm, respectively.  

\noindent Unlike most optimization algorithms using only the feasible solutions as the search domain, the FRI-FC problem benefits from infeasible points as well as the feasible ones. For this reason, the relaxed or perturbed inequalities (denoted by “\(\circ\)”) replace the ordinary ones (denoted by “\(\leq\)”) in the constraints of problem (1). Obviously, it is not reasonable to consider all the infeasible points with the same evaluation. For example, an infeasible point (however, with a good objective value) may be so far from the feasible region. Consequently, such a point cannot be reasonably considered as a solution for the problem. Thus, to evaluate the feasibility (or infeasibility) magnitude of the points, the feasibility property is defined as a fuzzy set so that closer points to the feasible set attain larger membership values. Figure 1(a) illustrates the notion of the fuzzy set “feasibility” by using grayscale colors varying from white to black in which the brighter parts include the points with larger membership values for “feasibility” (the white square – including all points between zero vector \(0\) and maximum solution \(\bar{x}\) –  depicts the feasible region of problem (2)). The consideration of such infeasible points as admissible selections for decision makers can be also viewed as perturbing some constraint(s) such that the feasible solutions set becomes greater and includes new points with likely better objective function values. In practice, the applied aspects of the problems sometimes allow some restrictions or limitations to be perturbed while preserving the purposes of the decision makers.

\noindent Similarly, the quality level or satisfaction amount of the objective value is defined as a fuzzy set. To define the fuzzy set “satisfaction”, the optimal objective value of problem (2), say \(z^{*}=c^{T}x^{*}\), is assumed to be the initial or current optimal value. Also, we consider a pre-determined value \(z_{0}(z_{0} < z^{*})\) as the best case and set its degree of membership (equivalently, its satisfaction amount) equal to one. Actually, the best case \(z_{0}\) is an approximate value we wish to be achieved at some near-feasible point (since for each feasible solution \(x\) in problem (2), we have \(z_{0} < z^{*} \leq z = c^{T}x\)). The value \(z_{0}\) may be determined f
rom some previous experiments or by the human experts. Anyway, the satisfaction amount for each objective value \(z = c^{T}x\) (especially, for \(z^{*}\) is realized by considering the difference between \(z_{0}\) and \(z\) (see Figure 1(b)). Ultimately, total values are obtained by evaluating the feasibility amount of the points as well as the quality level of their objective values and a point with the highest total value is introduced as an optimal solution. Briefly, the target of the FRI-FC problem is to find an infeasible point with the acceptable infeasibility and better objective value than \(z^{*}\) by taking advantage of the flexibility of constraint(s). Similar to [12,14], we refer to such these solutions as Super-Optima in this paper. If there is not a super-optimum, then the algorithm gives the same optimal solution obtained by other algorithms using only the feasible region as the search domain. In this case, in order to find a super-optimum, a decision maker has to consider more freedom for constraints to be perturbed (if it is possible based on the problem structure and the view of the human experts).

\begin{figure}[ht]
  \subcaptionbox*{(a)}[.40\linewidth]{%
    \includegraphics[height=\linewidth]{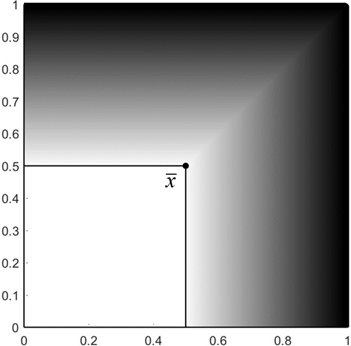}%
  }%
  \hfill
  \subcaptionbox*{(b)}[.40\linewidth]{%
    \hspace*{-2cm}\includegraphics[height=\linewidth]{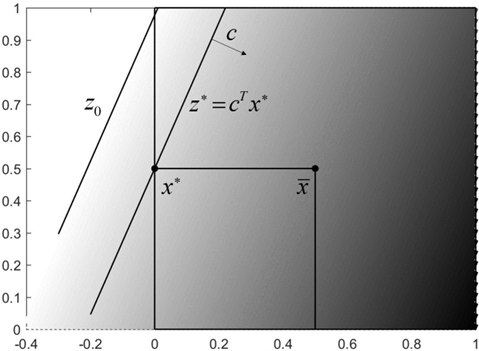}%
  }
  \caption{(a) Fuzzy evaluation of feasibility. (b) Fuzzy evaluation of objective values.}
\end{figure}

\noindent In this paper, we present a linearization approach for solving problem (1) in which \(\varphi\) is the product t-norm. In the proposed approach, the primal problem is initially converted into an equivalent linear optimization. By taking advantage of the linearity of the equivalent model, we can use many efficient methods such as Dantzig’s simplex algorithm, Karmarkar’s algorithm and interior-point methods for solving the problem. It is shown that the proposed linearization approach is an efficient and fast method that can find better solutions than those obtained by other related algorithms. On the contrary to the max-min algorithm [12] which stops once the first super-optimum is found, the proposed algorithm can find the best super-optimum. Moreover, in contrast to the modified PSO [14], the current approach provides exact solutions for the problem.

\noindent The rest of the paper is organized as follows. Section 2 takes a brief look at some basic results on the feasible solutions set of problem (2). In Section 3, the fuzzy constraints of problem (1) are precisely define by employing some membership functions. These membership functions are used for evaluating the feasibility and optimality amounts of points. Subsequently, problem (1) is transformed into an equivalent problem in which the constraints are linear. Two simplification rules are presented in Section 4 and in Section 5, the experimental results are demonstrated.

\noindent 
\vskip 0.2in
\noindent \textbf{2. Some previous theoretical results}
\vskip 0.2in

\noindent Consider problem (2) in which \(\varphi\) is the product t-norm and let \(x^{*}\) be its optimal solution with objective value \(z^{*}=c^{T}x^{*}\). Additionally, suppose there exist at least a flexible constraint in the problem, i.e., there are some points that violate the constraints and are still permissible to be solutions for the problem based on the decision maker’s view. Based on these assumptions, we can define problem (1) and find a better solution than \(x^{*}\) for problem (2) by taking advantage of the flexibility of some constraints. However, if the flexibility of the constraints is not sufficient, \(x^{*}\) is given again as an optimal solution. In this case, in order to find a super-optimum, a decision maker has to consider more freedom for constraints to be perturbed. At first, we mention some previously obtained results that are used throughout the paper. For the sake of simplicity, we let \(S(A, B)\) denote the feasible solutions set of problem (2), that is, \(S(A, B)=\{x \in [0,1]^{n}:A \varphi x \leq b\}\). So, similar to problem (1), we can rewrite \(S(A, B)\) in terms of rows \(a_{i}(i \in I)\) of matrix \(A\) as \(S(A, B)=\{x \in [0,1]^{n}:A \varphi x \leq b, i \in I\}\) where the constraints mean  \(a_{i} \varphi x = max_{j=1}^{n}\{\varphi (a_{ij}, x_{j})\} = max_{j=1}^{n}\{\varphi (a_{ij}, x_{i})\} \leq b_{i}, \forall i \in I\). 

\vskip 0.2in
\noindent \textbf{Definition 1.} \textbf{(a)} For each \(j \in J\) let \(I(j) = \{i \in I: a_{ij} > b_{j}\}\). \textbf{(b)} Let \(\bar{x} = ({\bar{x}}_{j})_{1 \times n}\) be an  dimensional vector in \(\) whose components are defined as \({\bar{x}}_{j} = \underset{i \in I (j)}{min}\{(b_{i}/a_{ij})\}, \forall j \in J\).

\vskip 0.2in

\noindent \textbf{Theorem 1.} \(S(A, B) = [0, \bar{x}]\). In other words, \(S(A, B)\) is a cube including all points \(x \in [0, 1]^{n}\) such that \(0 \leq x_{j} \leq \bar{x}_{j}, \forall j \in J\). (b) Solution \(x^{*} = ({x_1}^{*}, {x_2}^{*}, ..., {x_n}^{*})\) as defined below is the optimal solution for problem (2).

\begin{equation} \label{GrindEQ__1_} 
{x_{j}}^{*} =\begin{cases} \bar{x}_{j} & , c_{j} < b \\ 0 & , c_{j} \geq b\end{cases} 
\end{equation}

\noindent \textbf{Proof.} See Theorems 3 and 4 in [14].

\vskip 0.2in

\noindent \textbf{Theorem 2} (first simplification)\textbf{.} Consider problem (2) and let \(A' = (a'_{ij})_{m \times n}\) be a matrix resulted from matrix \(A = (a_{ij})_{m \times n}\) a follows. For \(i \in I\) and \(j_{0} \in J\), if there exist \(k \in I\) such that \(a_{i{j_{0}}} \geq b_{i}\),  \(a_{k{j_{0}}} > b_{k}\) and \(b_{i}/a_{i{j_{0}}} > b_{k}/a_{k{j_{0}}}\) , we set \(a'_{i{j_{0}}} = 0\); otherwise, \(a'_{i{j_{0}}} = a_{i{j_{0}}}\). Then, \(S(A', B) = S(A, B)\). 
\vskip 0.2in

\noindent \textbf{Proof.} Suppose that \(a_{i{j_{0}}} \geq b_{i}\),  \(a_{k{j_{0}}} > b_{k}\) and \(b_{i}/a_{i{j_{0}}} > b_{k}/a_{k{j_{0}}}\). We show that ‘‘resetting \(a_{i{j_{0}}}\) to zero’’ has no effect on \(S(A, B)\); that is, \(a_{i} \varphi x = max_{j=1}^{n}\{a_{ij} x_{j}\} \leq b_{i}\) is equivalent to \(\underset{j=1, j \neq j_{0}}{\overset{n}{max}}\{a_{ij} x_{j}\} < b_{i}, \forall x \in S(A, B)\). To this end, it is sufficient to prove \(a_{i{j_{0}}} x_{j_{0}} < b_{i}, \forall x \in S(A, B)\) . From \(max_{j=1}^{n}\{a_{ij} x_{j}\} \leq b_{i}\), it is obvious that \(a_{i{j_{0}}} x_{j_{0}} > b_{i}\) never holds. So, assume that \(a_{i{j_{0}}} x_{j_{0}} = b_{i}\). So, \(x_{j_{0}} = b_{i}/a_{i{j_{0}}} > b_{k}/a_{k{j_{0}}}\). Therefore, \(a_{k{j_{0}}} x_{j_{0}} > b_{k}\) which implies \(max_{j=1}^{n}\{a_{kj} x_{j}\} > b_{k}\). It contradicts \(x \in S(A, B)\) and the proof is complete.

\noindent 
\vskip 0.2in
\noindent \textbf{3. A mathematical model for Fuzzy constraints and equivalent problems}
\vskip 0.2in

\noindent In this section, we consider the constraints as well as the objective function of problem (1) as fuzzy sets by defining their associated membership functions. Via this approach, all points (whether feasible or not) in problem (2) are treated as feasible ones in (1) with different degrees of membership in the interval \([0,1]\). Actually, a point being in the feasible set of (2) belongs to that of (1) with the degree of membership equal to one, and a farther point from the feasible set of (2) gets a lower degree of membership. To treat the \(i\)'th fuzzy inequality of problem (1), we employ the same membership functions used in [12] as follows:

\begin{equation} \label{GrindEQ__1_} 
\mu(a_{i} \varphi x) =\begin{cases}
    1 & , a_{i} \varphi x \leq b_{i} \\
    1 - \frac{a_{i} \varphi x - b_{i}}{d_{i}} & , b_{i} \leq a_{i} \varphi x \leq b_{i} + d_{i} \\
    0 & , a_{i} \varphi x \geq b_{i} + d_{i}
\end{cases} \indent , i=1,2,...m
\end{equation}

\noindent where each \(d_{i}\) (\(i=1,2,...,m\) ) is initially determined as the limit of the admissible violation of the \(i\)‘th inequality. From relation (4), \(i\)‘th membership function is equal to 1 if the \(i\)‘th constraint of problem (2) is well satisfied, \(0\) if the constraint is violated beyond its admissible limit \(d_{i}\), and decreasing linearly from \(1\) to \(0\).

\vskip 0.2in
\noindent \textbf{Definition 2.} For each \(x \in [0,1]^{n}\) the crisp constraints violation vector at point \(x\) is an \(m \times 1\) vector \(CCV(x)=(CCV(x)_{i})_{m \times 1}\) whose \(i\)‘th component is defined as follows:

\begin{equation} \label{GrindEQ__1_} 
CCV(x)_{i} = max\{0, a_{i} \varphi x - b_{i}\}
\end{equation}

\noindent Also, we define the fuzzy constraints violation vector at point \(x\) as an \(m \times 1\) vector \(FCV(x)=(FCV(x)_{i})_{m \times 1}\) whose  \(i\)‘th component is

\begin{equation} \label{GrindEQ__1_} 
FCV(x)_{i} = max\{0, a_{i} \varphi x - (b_{i} + d_{i})\}
\end{equation}

\noindent From the above-mentioned definition, relations (5) and (6) show how much point \(x\) violates \(i\)‘th constraint of problem (2) (i.e., \(a_{i} \varphi x \leq b_{i}\)) and that of problem (1) (i.e., \(a_{i} \varphi x \circ b_{i}\)), respectively. Similar to (4), for the objective function of problem (1), we define

\begin{equation} \label{GrindEQ__1_} 
\mu(c^{T}x) =\begin{cases}
    1 & , c^{T}x \leq z_{0} \\
    1 - \frac{c^{T}x - z_{0}}{d_{0}} & , z_{0} \leq c^{T}x \leq z_{0} + d_{0} \\
    0 & , c^{T}x \geq z_{0} + d_{0}
\end{cases}
\end{equation}

\noindent where \(z_{0}=c^{T}x^{*}-vd_{0}\) for parameters \(v \in (0,1)\) and \(d_{0} \geq 0\) . From relation (7), we have \(\mu(z_{0})=1\), \(\mu(c^{T}x^{*})=1 - v\) and  \(\mu(z_{0} + d_{0})=0\). The parameters \(v\) and \(d_{0}\) determine the symmetry and length of the interval \([z_{0},z_{0}+d_{0}]=[c^{T}x^{*}-vd_{0},c^{T}x^{*}+(1-v)vd_{0}]\), respectively. If \(v \approx 0\) (\(v \approx 1\)) then \(\mu(c^{T}x^{*}) \approx 1\) (\(\mu(c^{T}x^{*}) \approx 0\)) that means the optimal solution of problem (2), \(x^{*}\), has a high (low) degree of satisfaction. If \(v=0.5\) then this interval is converted into the closure of the symmetric \(\varepsilon\)-neighborhood around \(c^{T}x^{*}\) with radius \(\varepsilon=0.5\). In this case, point \(x^{*}\) is interpreted as a solution that the satisfaction amount of its objective value is mediocre. Also, parameter \(d_{0}\) determines the length of the interval in which we want to survey attainability or non-attainability of a better objective function value than \(c^{T}x^{*}\); that is, we want to find a solution with a better objective function value than \(c^{T}x^{*}\) in this interval. Obviously, if \(d_{0}\) is selected very small (\(d_{0} \approx 0\)), then \(x^{*}\) can be considered as the best solution. In other words, this interval is actually an approximate domain in which we guess to be able to find a better objective function value than \(c^{T}x^{*}\) by partly perturbing some constraints. In [14], a discussion was provided on the initial setup of the parameters \(v\) and \(d_{0}\), and their influences on the convergence of the algorithm. Also, we refer the reader to [12] for a more detailed analysis of the influences of these parameters and some theoretical and experimental aspects that should be considered in problems.

\noindent By considering relations (4) and (7), we can express both the feasibility and optimality amounts of a point by only one variable as defined in the following definition.

\vskip 0.2in
\noindent \textbf{Definition 3.} For each \(x \in [0,1]^{n}\), let \(\mu_{0}(x)=\mu(c^{T}x)\) (i.e., the optimality value of \(x\)) and \(\mu_{F}(X)=min_{i=1}^{n}\{\mu(a_{i} \varphi x)\}\) (i.e., the feasibility value of \(x\)). Then, the total value of \(x\) is defined by \(\mu_{T}(X)=min\{\mu_{0}(x),mu_{F}(x)\}=min\{\mu(c^{T}x),\underset{i=1}{\overset{m}{min}}\{\mu(a_{i} \varphi x)\}\}\).

\noindent Definition 3 defines the whole space \(S=[0,1]^{n}\) as a fuzzy set by associating the membership value \(\mu(c^{T}x)\) to each \(x \in S\). Based on this fact, our purpose becomes equivalent to finding a point with the greatest total value among all \(x \in [0,1]^{n}\). In other words, we can express problem (1) as the following problem:

\begin{equation} \label{GrindEQ__1_} 
\underset{x \in S}{max}\{\mu_{T}(x)\}
\end{equation}

\noindent or equivalently

\begin{equation} \label{GrindEQ__1_} 
\underset{x \in [0,1]^{n}}{max}\{min\{\mu(c^{T}x),\underset{i=1}{\overset{m}{min}}\{\mu(a_{i} \varphi x)\}\}\}
\end{equation}

\noindent Therefore, by considering relations (4) and (7), problem (9) is rewritten as

\begin{equation} \label{GrindEQ__1_} 
\underset{x \in [0,1]^{n}}{max}\{min\{B_{0} - D_{0}(c^{T}x),\underset{i=1}{\overset{m}{min}}\{B_{0} - D_{0}(a_{i} \varphi x)\}\}\}
\end{equation}

\noindent where \(D_{i}=\frac{1}{d_{i}}\) for \(i \in I \cup {0}\), \(B_{0}=1+\frac{z_{0}}{d_{0}}\) and \(B_{i}=1+\frac{b_{i}}{d_{i}}\) for \(i \in I\). Now, by introducing the auxiliary variable \(\lambda\), problem (10) can be transformed into the following equivalent programming problem:

\begin{equation} \label{GrindEQ__1_} 
\begin{array}{l} {\max \, \, \, \, \lambda\, } \\ {\, \, \, \, \, \, \, \, \, \, \, \, D_{i}(a_{i} \varphi x) + \lambda \leq B_{i}\, ,i \in I} \\ {\, \, \, \, \, \, \, \, \, \, \, D_{0}(c^{T}x) + \lambda \leq B_{0}} \\ {\, \, \, \, \, \, \, \, \, \, \, x\in [0,1]^{n}} \end{array}
\end{equation}

\noindent or equivalently

\begin{equation} \label{GrindEQ__1_} 
\begin{array}{l} {\max \, \, \, \, \lambda\, } \\ {\, \, \, \, \, \, \, \, \, \, \, \, D_{i}(max_{j=1}^{n}\{a_{ij}x_{j}\}) + \lambda \leq B_{i}\, ,i \in I} \\ {\, \, \, \, \, \, \, \, \, \, \, D_{0}(c^{T}x) + \lambda \leq B_{0}} \\ {\, \, \, \, \, \, \, \, \, \, \, x\in [0,1]^{n}} \end{array}
\end{equation}

\noindent In problem (12), each constraint \(D_{i}(max_{j=1}^{n}\{a_{ij}x_{j}\})+\lambda \leq B_{i}\) is equivalent to \(n\) constraints \(D_{i}a_{ij}x_{j}+\lambda \leq B_{i}\) (\(j \in J\)). Due to this fact, problem (12) can be converted into the following optimization problem in which all the constraints are linear.

\begin{equation} \label{GrindEQ__1_} 
\begin{array}{l} {\max \, \, \, \, \lambda\, } \\ {\, \, \, \, \, \, \, \, \, \, \, D_{i}a_{ij}x_{j} + \lambda \leq B_{i}\, ,i \in I\, ,j \in J} \\ {\, \, \, \, \, \, \, \, \, \, \, D_{0}(c^{T}x) + \lambda \leq B_{0}} \\ {\, \, \, \, \, \, \, \, \, \, \, x\in [0,1]^{n}} \end{array}
\end{equation}

\noindent Heretofore, we showed that problem (1) is equivalent to problem (13) which is a linear programming problem. So, there are many efficient approaches such as Dantzig's simplex algorithm, Karmarkar’s algorithm and interior-point methods that can be used for solving problem (13).

\noindent As mentioned before, the modified PSO algorithm [14] was presented for solving problem (1) defined by an arbitrary continuous t-norm. Furthermore, by converting problem (1) into equivalent problem (8), we can use many heuristic algorithms which have been proposed for solving unconstrained optimization problems or those problems with some lower and upper bounds on the variables as their constraints. In Section 5, we apply the simplex algorithm to problem (13) and compare the generated results to those obtained by the modified PSO [14] and some well-known meta-heuristic methods which have been applied to many practical optimization problems.

\vskip 0.2in
\noindent \textbf{4. Simplification Rules}
\vskip 0.2in

\noindent As mentioned above, problems (12) and (13) are mathematically equivalent. So, by simplifying problem (12), we will have a more simplified linear optimization model (problem (13)) for solving. In this section, two simplification rules are presented to convert problem (12) into a more simplified equivalent form.

\vskip 0.2in
\noindent\textbf{Definition 4.} Let \(\lambda_{0}(x) = B_{0} - D_{0}(c^{T}x)\), \(\lambda_{i}(x) = B_{i} - D_{i}(max_{j=1}^{n}{a_{ij}x_{j}})\), \(\forall i \in I\), \(\lambda_{ij}(x_{j}) = B_{i} - D_{I}a_{ij}x_{j}\), \(\forall i \in I\)  and \(\forall j \in J\), and \(\Lambda(x) = max_{i=0}^{m}{\lambda_{i}(x)}\).

\vskip 0.2in
\noindent \textbf{Corollary 1.} Functions \(\lambda_{i}(x)\), \(\forall i \in I\), are non-increasing continuous functions. Functions \(\lambda_{ij}(x_{j})\), \(\forall i \in I\) and \(\forall j \in J\), are strictly decreasing and function \(\lambda_{0}(x)\) is strictly decreasing with respect to component \(x_{j}\) if \(c_{j} > 0\) and strictly increasing if \(c_{j} < 0\).

\noindent By Definition 4, problem (12) can be expressed as follows:

\begin{equation} \label{GrindEQ__1_} 
\underset{x \in [0,1]^{n}}{max}\{\underset{i=0}{\overset{m}{min}}\{\lambda_{i}(x)\}\}
\end{equation}

\noindent or equivalently,

\begin{equation} \label{GrindEQ__1_} 
\underset{x \in [0,1]^{n}}{max}\{\Lambda(x)\}
\end{equation}

\noindent Theorem 3 provides a simplification rule for solving problem (14) by finding some components of the optimal solution.

\vskip 0.2in
\noindent\textbf{Theorem 3.} Suppose that \(x^{*}\) is an optimal solution of problem (14). Then, \({x_{j}}^{*} = 0\) for each \(j \in J\) such that \(c_{j} > 0\).

\vskip 0.2in
\noindent \textbf{Proof.} Since (12) and (15) are equivalent problems, \(x^{*}\) is also an optimal for (15). Thus, \(\Lambda(x^{*}) = \Lambda(x)\) for each \(x \in [0,1]^{n}\). By contradiction, suppose that \({x_{j_{0}}}^{*} > 0\) for some \(j_{0} \in J\) such that \(c_{j_{0}}^{*} > 0\). Let \(x' \in [0,1]^{n}\) such that \(x'_{j_{0}} = {x_{j_{0}}}^{*}\) if \(c_{j} < 0\) and \(x'_{j_{0}} = 0\) if \(c_{j} > 0\). Then, from Corollary 1, we have \(\lambda_{i}(x') \geq \lambda_{i}(x^{*})\), \(\forall i \in I \cap \{0\}\) , and therefore \(\Lambda(x') = \Lambda(x^{*})\) that contradicts the optimality of \(x^{*}\).

\vskip 0.2in
\noindent\textbf{Corollary 2.} (second simplification). For each \(j \in J\) such that \(c_{j} > 0\), we can assign \({x_{j}}^{*} = 0\) and remove the \(j\)‘th column of matrix \(A\). Therefore, in order to solve problem (12) (or problem (13)), it is sufficient to consider only the columns of matrix \(A\) that belong to \(J' = \{j \in J : c_{j} < 0\}\).

\vskip 0.2in
\noindent\textbf{Lemma 1.} Let \(i \in I\). Then, \(\lambda_{i}(x) = min_{j in J}\{\lambda_{ij}(x_{j})\}\), \(\forall x \in [0,1]^{n}\) . 

\vskip 0.2in
\noindent \textbf{Proof.} Let \(x \in [0,1]^{n}\) and \(\lambda_{ij_{0}}(x_{j_{0}}) = min_{j in J}\{\lambda_{ij}(x_{j})\}\). So, from Definition 4 we have \(B_{i} - D_{i}a_{ij_{0}}x_{j_{0}} \leq B_{i} - D_{i}a_{ij}x_{j}\), \(\forall j \in J\) . Therefore, \(a_{ij}x_{j} \leq a_{ij_{0}}x_{j_{0}}\),  \(\forall j \in J\), which together with Definition 4 imply that \(\lambda_{i}(x) = B_{i} - D_{i}(max_{j=1}^{n}{a_{ij}x_{j}}) = B_{i} - D_{i}a_{ij_{0}}x_{j_{0}} = \lambda_{ij_{0}}(x_{j_{0}})\) .

\noindent From Corollary 2 and Lemma 1, problem (14) is converted into the following simplified problem:

\begin{equation} \label{GrindEQ__1_} 
\underset{x \in [0,1]^{n}}{max}\{min\{\lambda_{0}(x), \underset{i \in I, j \in J}{min}\{\lambda_{ij}(x_{j})\}\}\}
\end{equation}


\vskip 0.2in
\noindent\textbf{Lemma 2.} Let \(i \in I\) and \(j \in J'\). If \(a_{ij} \leq b_{i}\) , then \(\lambda_{ij}(x_{j}) \geq 1\) , \(\forall x_{j} \in [0,1]^{n}\).

\vskip 0.2in
\noindent \textbf{Proof.} Suppose that \(a_{ij} \geq b_{i}\). Thus, \(a_{ij}x_{j} \geq b_{i}\) and therefore \(1 - \frac{a_{ij} - b_{i}}{d_{i}} \geq 1\), \(\forall x_{j} \in [0,1]\) . Now, by substituting \(D_{i} = \frac{1}{d_{i}}\) and \(B_{i} = 1 + \frac{b_{i}}{d_{i}}\) in Definition 4, we have \(\lambda_{ij}(x_{j}) = 1 - \frac{a_{ij} - b_{i}}{d_{i}} \geq 1\). 

\vskip 0.2in
\noindent\textbf{Theorem 4.} Suppose the second simplification (Corollary 2) is done and \(i \in I\). Then, \(\lambda_{i}(x) = min_{j \in J'_{i}}\{\lambda_{ij}(x_{j})\}\) where \(J'_{i} = \{j \in J' : a_{ij} > b_{i}\}\) .

\vskip 0.2in
\noindent \textbf{Proof.} From Corollary 2 and Lemma 1, \(\lambda_{i}(x) = min_{j \in J}\{\lambda_{ij}(x_{j})\}\), \(\forall i \in I\). Hence, \(\lambda_{i}(x) = min\{min_{j \in J'_{i}}\{\lambda_{ij}(x_{j})\}, min_{j \notin J'_{i}}\{\lambda_{ij}(x_{j})\}\}\) , \(\forall i \in I\). Since \(j \notin J'_{i}\) implies \(a_{ij} \leq b_{i}\), from Lemma 2 we obtain \(\lambda_{ij}(x_{j}) \geq 1\), \(\forall j \in J' - j'_{i}\) and \(\forall x \in [0,1]\). Therefore, \(min_{j \in J'_{i}}\{\lambda_{ij}(x_{j})\} \geq 1\) and then \(\lambda_{i}(x) = min_{j \in J'_{i}}\{\lambda_{ij}(x_{j})\}\).

\noindent From Theorem 4, problem (16) is converted into the more simplified form (17) as follows:

\begin{equation} \label{GrindEQ__1_} 
\underset{x \in [0,1]^{n}}{max}\{min\{\lambda_{0}(x), \underset{i \in I, j \in J}{min}\{\lambda_{ij}(x_{j})\}\}\}
\end{equation}

\vskip 0.2in
\noindent\textbf{Corollary 3} (third simplification) \textbf{. (a)} Suppose that \(i_{0} \in I\) and \(j_{0} \notin J'_{i}\). Then, ‘’resetting \(a_{i_{0}j_{0}}\) to zero’’ has no effect on the optimal solution of problem (14) (or that of problem (12). As a result, if \(j_{0} \notin J'_{i}\), \(\forall i \in I\), then we can remove \(j_{0}\)‘th column of matrix \(A\). \textbf{(b)} Let \(x^{*}\) be an optimal solution of problem (14). \(j_{0} \notin J'_{i}\), \(\forall i \in I\), then \({x_{j_{0}}}^{*} = 1\).

\vskip 0.2in
\noindent \textbf{Proof.} \textbf{a} Since \(j_{0} \notin J'_{i_{0}}\), then \(\lambda_{i_{0}j_{0}}(x_{j_{0}})\) does not appear in (17) and therefore we can take \(\lambda_{i_{0}j_{0}}(x_{j_{0}})\) out of consideration. So, the value of \(a_{i_{0}j_{0}}\) plays no role in (17) and can be assigned to zero. \textbf{b} By contradiction, suppose \({x_{j_{0}}}^{*} < 1\) . Let \(x' \in [0,1]^{n}\) such that \(x'_{j_{0}} = 1\) and \(x'_{j} = {x_{j}}^{*}\), \(\forall j \in J'_{i}\) . Then, \(\underset{in \in I, j \in J'_{i}}{min}\{\lambda_{ij(x'_{j})}\} = \underset{in \in I, j \in J'_{i}}{min}\{\lambda_{ij({x_{j}}^{*})}\}\) and we also have \(\lambda_{i}(x') \geq \lambda_{i}(x^{*})\) by Corollary 1. Hence, \(\Lambda(x') = \Lambda(x^{*})\) that contradicts the optimality of \(x^{*}\) .

\vskip 0.2in
\noindent \textbf{5. Comparisons with other works and some numerical examples}
\vskip 0.2in

\noindent In this section, we present the experimental results for evaluating the performance of our algorithm. In Section 5.1, the proposed algorithm is applied for solving the transformed equivalent problem (problem (13)) for the test problems described in Appendix A.  The test problems have been randomly generated in different sizes by using product t-norm. For each test problem, we have a pair of problems that are associated with each other; one FRI-FC problem (problem (1)) and one FRI problem with ordinary (crisp) inequalities (problem (2)). From Theorem 1(a), we know that problem (2) is feasible and \(S(A, B) = [0, \bar{x}]\) where the maximum solution \(\bar{x}\) is obtained from Definition 1. Using this result, the feasible optimal solution \(x^{*}\) of problem (2) is directly given by Theorem 1(b). To perform a fair comparison, we use the same parameters for each experiment, that is, we set \(v = 0.5\) and \(d_{i} = 0.1\) for \(i = 0, 1,..., m\). Therefore, from (4) and (7) we have \(\mu_{T}(z_{0}) = 1\) and \(\mu_{T}(c^{T}x^{*}) = 1 - v = 0.5\). The equality \(\mu_{T}(c^{T}x^{*}) = 0.5\) means that \(x^{*}\) (the optimal solution of problem (2)) has a mediocre objective value for us. Actually, the target of the methods is to find a solution \(x^{**}\) such that \(c^{T}x^{**} \in [z_{0}, z_{0} + d_{0}]\) and \(c^{T}x^{**}\) is as close as possible to \(z_{0}\) (the best case) by perturbing the constraints within the range determined by \(d_{i} = 0.1\), \(i = 0, 1,..., m\). Finally, a comparison is made between the current method, the Modified PSO (MPSO) [14], Original PSO (OPSO) [24], Continuous Ant Colony Optimization (CACO) [45], Differential Evolution (DE) [38] and Harmony Search (HS) [11] algorithms. For this purpose, in Section 5.2, we apply these methods to the ten test problems described in Appendix A.

\noindent \\ 5.1. Results of the linearization approach

\noindent In this section, the proposed algorithm is applied to the test problems described in Appendix A. Table 1 includes the feasibility values \(\mu_{F}(x^{**})\), optimality values \(\mu_{0}(x^{**})\) and total values \(\mu_{T}(x^{**})\) for the best solutions \(x^{**}\) found by the linearization algorithm. In each case, to simplify the comparison between the current optimum \(x^{*}\) (the optimal solution for problem (2)) and the best super-optimum \(x^{**}\) found by the proposed algorithm, values \(c^{T}x^{*}\) and \(c^{T}x^{**}\) have been reported. Also, we reported admissible ranges \([z_{0}, z_{0} + \varepsilon_{0}]\) for each test problem to determine if the proposed method can find a super-optimum \(x^{**} \in [z_{0}, z_{0} + d_{0}]\) such that \(c^{T}x^{**} < c^{T}x^{*}\) . Additionally, for each test problem, Table 2 presents the optimal solutions \(x^{*}\) of problem (2), solutions \(x^{**}\) found by the current algorithm and the crisp constraints violation vectors at points  \(x^{**}\). \\

\begin{table}
\caption{\label{demo-table}A comparison between the quality of the best super-optima (i.e., solutions \(x^{**}\)) found by the proposed linearization approach and that of the optimal solutions of problem (2) (i.e., solutions \(x^{*}\)) for the test problems A.1 - A.10.}
\begin{center}
\resizebox{\columnwidth}{!}{
\begin{tabular}{|c|c|c|c|c|c|c|} 
 \hline
Test problems  & \(\mu_{F}(x^{**})\) & \(\mu_{0}(x^{**})\) & \(\mu_{T}(x^{**})\) & \(c^{T}x^{**}\) & \(c^{T}x^{*}\) & \([z_{0}, z_{0} + d_{0}]\) \\ \hline
A.1 & 0.9910 & 0.9910 & 0.9910 & - 0.9232 & - 0.8741 & [- 0.9241 , - 0.8241] \\
\hline
A.2 & 0.9933 & 0.9933 & 0.9933 & - 11.3722 & - 11.3228 & [- 11.3728 , - 11.2728] \\
\hline
A.3 & 0.9765 & 0.9765 & 0.9765 & - 1.3501 & - 1.3024 & [- 1.3524 , - 1.2524] \\
\hline
A.4 & 0.9916 & 0.9916 & 0.9916 & - 9.7886 & - 9.7395 & [- 9.7895 , - 9.6895] \\
\hline
A.5 & 0.9793 & 0.9793 & 0.9793 & - 1.4395 & - 1.3916 & [- 1.4416 , - 1.3416] \\
\hline
A.6 & 0.9809 & 0.9809 & 0.9809 & - 0.1638 & - 0.1157 & [- 0.1657 , - 0.0657] \\
\hline
A.7 & 0.9647 & 0.9647 & 0.9647 & - 0.0820 & - 0.0356 & [- 0.0856 , 0.0144] \\
\hline
A.8 & 0.9373 & 0.9373 & 0.9373 & - 0.1337 & - 0.0899 & [- 0.14 , - 0.04] \\
\hline
A.9 & 0.9907 & 0.9907 & 0.9907 & - 1.0552 & - 1.0061 & [- 1.0561 , - 0.9561] \\
\hline
A.10 & 0.9914 & 0.9914 & 0.9914 & - 1.7222 & - 1.6731 & [- 1.7231 , - 1.6231] \\
 \hline
\end{tabular}
}
\end{center}
\end{table}

\noindent As shown in Table 1, we have \(\mu_{T}(x^{**}) > \mu_{T}(c^{T}x^{*}) = 0.5\) in all the cases,  that means the proposed algorithm can find solutions with higher quality than \(x^{*}\). According to Table 1, objective values \(c^{T}x^{**}\) belong to the admissible intervals \([z_{0}, z_{0} + d_{0}]\) (and very close to the best cases \(z_{0}\)) and are strictly less than \(c^{T}x^{*}\) for all the test problems. Therefore, the linearization method produces better solutions (i.e., solutions with less objective values) compared to the optimal solutions of problem (2). Table 2 includes solutions \(x^{*}\), \(x^{**}\) and the crisp constraints violation vectors \(CCV(x^{**})\) for each test problem. The results, in Table 2, show that the proposed method produces optimal solutions with admissible infeasibilities; more precisely, for each \(i \in I\), we have \(CCV(x^{**})_{i} < d_{i} = 0.1\). As a key result, we see that the best super-optimum \(x^{**}\) satisfies \(\mu_{F}(x^{**}) = \mu_{0}(x^{**})\) for each test problem. As it will be shown in the next section, this equality does not hold true for each super-optimum.

\begin{table}
\caption{\label{demo-table}The best super-optima, \(x^{**}\) , and their crisp constraints violation vectors \(CCV(x^{**})\).}
\begin{center}
\resizebox{\columnwidth}{!}{
\begin{tabular}{|c|l|} 
 \hline
A.1 & \(x^{**}=[0.1859,0.115,0.0165,0,0,0]\) \\  & \(x^{**}=[0.1964,0.1215,0.0174,0,0,0]\) \, \(CCV(x^{**})=[0,0.0009,0,0]\) \\
\hline
A.2 & \(x^{**}=[0,0.8702,0,0,0.0165,0.5835]\) \\  & \(x^{**}=[0,0.8731,0,0,0.6506,0.5854]\) \, \(CCV(x^{**})=[0,0.0007,0.0007,0,0]\) \\
\hline
A.3 & \(x^{**}=[0,0.0979,0,0.1334,0,0.1532]\) \\  & \(x^{**}=[0,0.1015,0,0.1383,0,0.1588]\) \, \(CCV(x^{**})=[0,0,0,0,0.0023,0]\) \\
\hline
A.4 & \(x^{**}=[0.2069,0,0,0.0339,0,0,0.8705]\) \\  & \(x^{**}=[0.2123,0,0,0.0348,0,0,0.8717]\) \, \(CCV(x^{**})=[0.0008,0,0,0.0008,0]\) \\
\hline
A.5 & \(x^{**}=[0.3282,0,0,0.1228,0,0,0.1221,0]\) \\  & \(x^{**}=[0.3319,0,0,0.1302,0,0,0.1294,0]\) \, \(CCV(x^{**})=[0,0.0021,0.0021,0,0,0]\) \\
\hline
A.6 &\(x^{**}=[0,0,0,0,0.0058,0.0067,0.0056]\) \\  & \(x^{**}=[0,0,0,0,0.0082,0.0094,0.0079]\) \, \(CCV(x^{**})=[0,0,0,0.0019,0,0,0]\) \\
\hline
A.7 & \(x^{**}=[0.0057,0,0,0,0,0,0,0,0,0.0043]\) \\  & \(x^{**}=[0.0131,0,0,0,0,0,0,0,0.0.0099]\) \, \(CCV(x^{**})=[0,0,0,0,0,0,0.0035]\) \\
\hline
A.8 & \(x^{**}=[0.0613,0,0,0,0,0.0163,0,0,0,0]\) \\  & \(x^{**}=[0.0910,0,0,0,0,0.0242,0,0,0,0]\) \, \(CCV(x^{**})=[0,0,0,0,0.0063,0,0,0]\) \\
\hline
A.9 &\(x^{**}=[0,0.0234,0.1036,0,0.0235,0.0228,0,0,0]\) \\  & \(x^{**}=[0,0.0245,0.1087,0,0.0.246,0.0240,0,0,0]\) \, \(CCV(x^{**})=[0,0.0009,0,0,0,0,0,0,0]\) \\
\hline
A.10 & \(x^{**}=[0,0.0295,0,0,0.0335,0.0633,0,0,0,0.1846]\) \\  & \(x^{**}=[0,0.0304,0,0,0.0345,0.0652,0,0,0,0.1900]\) \, \(CCV(x^{**})=[0.0009,0,0,0,0,0,0,0,0,0]\) \\
 \hline
\end{tabular}
}
\end{center}
\end{table}

\noindent \\ 5.2. Linearization approach versus the other related methods

\noindent In this section, a comparison is made between the current linearization method and the modified PSO algorithm [14] proposed for solving the FRI-FC problems. As mentioned before, since problem (1) is equivalent to problem (9), many heuristic algorithms may be used for solving the problem. So, the generated solutions for the linearization method and modified PSO are also compared with some well-known meta-heuristics such as Original PSO (OPSO) [24], Continuous Ant Colony Optimization (CACO) [45], Differential Evolution (DE) [38] and Harmony Search (HS) [11] algorithms. For the heuristic algorithms, \(30\) experiments are performed for each test problem. The maximum number of iterations is equal to \(100\). The parameters of the PSO algorithms that are used in each case are as follows. Swarm size is set to \(10\), \(c_{1} = C_{2} = 2\) and inertia factor \(w\) is decreasing linearly from \(1\) to \(0\) by a damping factor of  [14,24]. For CACO, \(m = 10\) (number of ants used in an iteration), \(\xi = 0.85\) (the speed of convergence), \(q = 10^{-4}\) (locality of the search process) and \(k = 50\) (archive size) [45], and for DE, \(c_{T} = 0.95\) (crossover rate) and \(F = 0.5\) (mutation scaling factor) [38], and for HS, \(HMS = 10\) (size of harmony memory), \(HMCR = 0.825\) (harmony memory consideration rate) and \(PAR = 0.35\) (Pitch Adjustment Rate) [11].

\noindent Table 3 presents the total values of the best solutions obtained by the algorithms. As shown in this table, CACO, DE and HS have the worst results (i.e., solutions \(x\) with the least total values \(\mu_{T}(x) = 0\)) in most cases (\(70\%\) of the test problems).  Moreover, for the remaining \(30\%\) of the cases, they could not necessarily find a super-optimum (i.e., solutions \(x\) such that \(\mu_{T}(x) > \mu_{T}(x^{*})\) ); A.3 for CACO, and A.4 for DE and HS. Also, OPSO produced the worst solutions in half of the cases and could find super-optima only for \(50\%\) of the test problems. On the other hand, the modified PSO could find a super-optimum for each test problem that is very close to the best super-optimum \(x^{**}\) obtained by the linearization approach. Particularly, for the test problem A.8, we have \(\mu_{T}(x^{*}_{MPSO}) = \mu_{T}(x^{**})\) . However, for test problem A.4, although \(\mu_{T}(x^{*}_{MPSO}) > \mu_{T}(x^{*})\), the distance between \(\mu_{T}(x^{*}_{MPSO})\) and \(\mu_{T}(x^{**})\) is very large in average.

\noindent In Table 4, the results have been averaged over 30 runs and the average best-so-far (Avg), median best-so-far (Mdn) in the last iterations and the standard deviations (Sd) are reported for the heuristic algorithms. The results, in Table 4, show the modified PSO produces better solutions with a higher convergence rate when compared against the other heuristic algorithms. As this table illustrates, \(70\%\) of the average best-so-far values found by CACO, DE and HS are zero. As mentioned earlier, this is because of the unfavorable points \(x\) having zero total values. Due to this fact, it can be concluded that these algorithms could not escape from these poor solutions for \(70\%\) of the problems. Also, OPSO could avoid trapping in such unfavorable points only for half of the cases. Figures 2-8 depict the performances of the algorithms that are averaged over 30 runs for some test problems. Dashed lines \(y = \mu_{T}(x^{*})\) and \(y = \mu_{T}(x^{**})\) determines the least and the most accuracy required, respectively.

\begin{table}
\caption{\label{demo-table}A comparison between the total values of the best solutions \(x^{**}\) , \(x^{*}_{MPSO}\) , \(x^{*}_{OPSO}\) , \(x^{*}_{CACO}\) , \(x^{*}_{DE}\)  and \(x^{*}_{HS}\)  found by the linearization approach, MPSO, OPSO, CACO, DE and HS algorithms, respectively, for the test problems of Appendix A.}
\begin{center}
\resizebox{\columnwidth}{!}{
\begin{tabular}{|c|c|c|c|c|c|c|} 
 \hline
Test problems  & \(\mu_{T}(x^{**})\) & \(\mu_{T}(x^{*}_{MPSO})\) & \(\mu_{T}(x^{*}_{OPSO})\) & \(\mu_{T}(x^{*}_{CACO})\) & \(\mu_{T}(x^{*}_{DE})\) & \(\mu_{T}(x^{*}_{HS})\) \\ \hline

A.1 & 0.9910 & 0.9904 & 0.9608 & 0.9632 & 0.5638 & 0.7519 \\
\hline
A.2 & 0.9933 & 0.9924 & 0.9929 & 0.9895 & 0.9480 & 0.9480 \\
\hline
A.3 & 0.9765 & 0.9717 & 0.6582 & 0.0834 & 0 & 0 \\
\hline
A.4 & 0.9916 & 0.9901 & 0.9494 & 0 & 0.2619 & 0.1217 \\
\hline
A.5 & 0.9793 & 0.9749 & 0 & 0 & 0 & 0 \\
\hline
A.6 & 0.9809 & 0.9788 & 0.9390 & 0 & 0 & 0 \\
\hline
A.7 & 0.96472 & 0.96471 & 0 & 0 & 0 & 0 \\
\hline
A.8 & 0.9373 & 0.9373 & 0 & 0 & 0 & 0 \\
\hline
A.9 & 0.9907 & 0.9870 & 0 & 0 & 0 & 0 \\
\hline
A.10 & 0.9914 & 0.9875 & 0 & 0 & 0 & 0 \\
 \hline
\end{tabular}
}
\end{center}
\end{table}

\noindent Table 5 contains the error of the solutions found by the algorithms (for heuristic algorithms, the results are the average error of 30 runs over the last iterations). And the error is given by \(\frac{1}{2}\{(\frac{1}{m}) \sum_{i=1}^m CCV(x_{i}) + |z_{0} - c^{T}x|\}\), where \((\frac{1}{m}) \sum_{i=1}^m CCV(x_{i})\) is the average of all the crisp constraint violations, and \(|z_{0} - c^{T}x|\) is the difference between \(z_{0}\) (the best desirable value) and optimal values obtained by applying the algorithms to the test problems. From this table, we can note that the linearization approach produces the best super-optima with the least errors and therefore with the least constraint violations (also, see Table 2). Additionally, the time (in seconds) involved in the execution of each method is reported in Tables 6. Based on the results, in Table 6, the current approach is the fastest algorithm in comparison to the other algorithms in each test problem. Actually, the linear approach takes less than 0.996 seconds to finish an experiment.

\begin{table}
\caption{\label{demo-table}A comparison between the results found by the MPSO, OPSO, CACO, DE and HS algorithms for the test problems of Appendix A. The results have been averaged over 30 runs. Maximum number of iterations=100.}
\begin{center}
\resizebox{\columnwidth}{!}{
\begin{tabular}{|c|c|c|c|c|c|c|} 
\hline
Test problems  &   & MPSO & OPSO & CACO & DE & HS \\ \hline
     & Avg & 0.9759 & 0.091549 & 0.059136 & 0.018792 & 0.045539  \\ 
     \cline{2-7}
 A.1 & Mdn & 0.9787 & 0 & 0 & 0 & 0 \\
     \cline{2-7}
     & Sd & 0.0035 & 2.7014 & 1.3873 & 1.3413 & 1.5425  \\ 
     \hline \hline
     & Avg & 0.93161 & 0.34901 & 0.6167 & 0.49248 & 0.49248  \\ 
     \cline{2-7}
 A.2 & Mdn & 0.97558 & 0 & 0.84622 & 0.57185 & 0.57185 \\
     \cline{2-7}
     & Sd & 0.58775 & 2.8787 & 2.5854 & 0.89195 & 0.89195  \\ 
     \hline \hline     
     & Avg & 0.83368 & 0.038722 & 0.0027787 & 0 & 0  \\ 
     \cline{2-7}
 A.3 & Mdn & 0.87382 & 0 & 0 & 0 & 0 \\
     \cline{2-7}
     & Sd & 0.62053 & 2.7763 & 1.885 & 1.9039 & 1.7758  \\ 
     \hline \hline
     & Avg & 0.5177 & 0.031648 & 0 & 0.0087283 & 0.0040561 \\ 
     \cline{2-7}
 A.4 & Mdn & 0.50478 & 0 & 0 & 0 & 0 \\
     \cline{2-7}
     & Sd & 0.40083 & 2.1511 & 1.8494 & 1.8822 & 2.1892  \\ 
     \hline \hline
     & Avg & 0.92709 & 0 & 0 & 0 & 0  \\ 
     \cline{2-7}
 A.5 & Mdn & 0.96284 & 0 & 0 & 0 & 0 \\
     \cline{2-7}
     & Sd & 0.092625 & 4.3403 & 2.0897 & 1.9291 & 2.5145  \\ 
     \hline \hline
     & Avg & 0.93838 & 0.031301 & 0 & 0 & 0  \\ 
     \cline{2-7}
 A.6 & Mdn & 0.95304 & 0 & 0 & 0 & 0 \\
     \cline{2-7}
     & Sd & 0.0024304 & 1.6968 & 1.4244 & 1.505 & 1.68  \\ 
     \hline \hline
     & Avg & 0.96435 & 0 & 0 & 0 & 0  \\ 
     \cline{2-7}
 A.7 & Mdn & 0.96422 & 0 & 0 & 0 & 0 \\
     \cline{2-7}
     & Sd & 1.2739e-05 & 4.0588 & 2.4427 & 2.0715 & 2.7796  \\ 
     \hline \hline
     & Avg & 0.93015 & 0 & 0 & 0 & 0  \\ 
     \cline{2-7}
 A.8 & Mdn & 0.93718 & 0 & 0 & 0 & 0 \\
     \cline{2-7}
     & Sd & 0.0020195 & 5.0862 & 3.454 & 3.3177 & 2.7028  \\ 
     \hline \hline
     & Avg & 0.94184 & 0 & 0 & 0 & 0  \\ 
     \cline{2-7}
 A.9 & Mdn & 0.92721 & 0 & 0 & 0 & 0 \\
     \cline{2-7}
     & Sd & 0.0027109 & 2.4319 & 1.4109 & 1.1276 & 0.90631  \\ 
     \hline \hline
     & Avg & 0.91627 & 0 & 0 & 0 & 0  \\ 
     \cline{2-7}
 A.10 & Mdn & 0.95904 & 0 & 0 & 0 & 0 \\
     \cline{2-7}
     & Sd & 0.0082584 & 2.3339 & 1.7099 & 1.796 & 1.788  \\ 
     \hline
     
\end{tabular}
}
\end{center}
\end{table}

\noindent Table 7 shows the results of paired t-test to illustrate significant difference between the algorithms. In the pair t-test, the result is a value called p-value, and if it is close to zero, it indicates significance difference between two methods. Also, confidential interval implies interval of difference in which  \(95\%\) of values lie. In Table 7, Sd, SEM and CI stand for standard deviation of differences, standard error mean of differences and confidential interval, respectively. As illustrated in Table 7, the p-values are mostly lower than \(10^{-7}\) which indicates significant difference between the algorithms. These results show that the differences between the optimal solutions found by the OPSO, CACO, DE and HS algorithms and the best super-optima found by the linearization approach are statistically significant. Moreover, for the modified PSO, the results indicate the closeness between these solutions and the best super-optima.

\begin{table}
\caption{\label{demo-table}Errors averaged over the last iterations of 30 runs.}
\begin{center}
\resizebox{\columnwidth}{!}{
\begin{tabular}{|c|c|c|c|c|c|c|} 
 \hline
Test problems  & Proposed Method & MPSO & OPSO & CACO & DE & HS \\ \hline

A.1 & 0.00056521 & 0.0022895 & 2.9873 & 2.1724 & 2.6526 & 2.3742 \\
\hline
A.2 & 0.00047085 & 0.12277 & 2.8629 & 1.4587 & 0.52682 & 0.52682 \\
\hline
A.3 & 0.0013702 & 0.12296 & 3.6869 & 3.3445 & 3.5706 & 3.9572 \\
\hline
A.4 & 0.00058783 & 0.31098 & 3.3696 & 3.0948 & 3.1164 & 3.6959 \\
\hline
A.5 & 0.0013822 & 0.019544 & 7.9953 & 7.3276 & 7.544 & 7.7492\\
\hline
A.6 & 0.0010926 & 0.0035921 & 2.3294 & 1.6659 & 2.3062 & 2.6638\\
\hline
A.7 & 0.0020163 & 0.0020349 & 8.5333 & 8.5898 & 8.7956 & 9.5084 \\
\hline
A.8 & 0.0035273 & 0.0038748 & 10.7297 & 13.4616 & 12.9967 & 12.9733 \\
\hline
A.9 & 0.00051485 & 0.003506 & 2.9832 & 2.0228 & 1.9556 & 1.7978 \\
\hline
A.10 & 0.00047489 & 0.0061162 & 3.6554 & 2.6794 & 2.7964 & 2.9454\\
 \hline
 \hline
MSE & 2.2866e-06 & 0.012737 & 32.1624 & 35.035 & 34.9926 & 37.3003\\
\hline
\end{tabular}
}
\end{center}
\end{table}

\begin{table}
\caption{\label{demo-table}Time involved in the execution of the methods.}
\begin{center}
\resizebox{\columnwidth}{!}{
\begin{tabular}{|c|c|c|c|c|c|c|} 
 \hline
Test problems  & Proposed Method & MPSO & OPSO & CACO & DE & HS \\ \hline

A.1 &  0.036123 & 0.28818 & 0.15882 & 0.19304 & 0.14393 & 0.20819\\
\hline
A.2 &  0.099543 & 0.31729 & 0.17306 & 0.2198 & 0.1677 & 0.1677\\
\hline
A.3 &  0.071229 & 0.33334 & 0.18697 & 0.23123 & 0.17436 & 0.2457\\
\hline
A.4 &  0.063235 & 0.35112 & 0.19826 & 0.23182 & 0.17636 & 0.25039\\
\hline
A.5 & 0.092431 & 0.3857 & 0.21215 & 0.2583 & 0.19873 & 0.28562\\
\hline
A.6 & 0.088895 & 0.39704 & 0.21742 & 0.26146 & 0.20715 & 0.28364\\
\hline
A.7 &  0.041282 & 0.55216 & 0.25899 & 0.31076 & 0.24577 & 0.34284\\
\hline
A.8 &  0.026400 & 0.516 & 0.28584 & 0.33663 & 0.27376 & 0.3705\\
\hline
A.9 & 0.042752 & 0.5221 & 0.28127 & 0.33833 & 0.27599 & 0.36654\\
\hline
A.10 & 0.035890 & 0.59621 & 0.32267 & 0.38421 & 0.32517 & 0.41314\\
 \hline
\end{tabular}
}
\end{center}
\end{table}

\begin{table}
\caption{\label{demo-table}Paired t-test results based on output value and with degree freedom of 7.}
\begin{center}
\resizebox{\columnwidth}{!}{
\begin{tabular}{|c|c|c|c|c|c|c|} 
 \hline
Test & Difference between & Difference between & Difference between & Difference between &  Difference between \\
 problems & proposed method & proposed method & proposed method &  proposed method & proposed method \\
& and MPSO &  and OPSO &  and CACO & and DE &  and HS \\ \hline

A.1 &  0.0151 & 0.89941 & 0.93182 & 0.97217 & 0.94542\\
\hline
A.2 &  0.06166 & 0.64426 & 0.37657 & 0.50079 & 0.50079\\
\hline
A.3 & 0.14283 & 0.93779 & 0.97373 & 0.97651 & 0.97651\\
\hline
A.4 & 0.4739 & 0.95995 & 0.9916 & 0.98287 & 0.98754\\
\hline
A.5 & 0.05218 & 0.97927 & 0.97927 & 0.97927 & 0.97927\\
\hline
A.6 & 0.0425 & 0.94958 & 0.98088 & 0.98088 & 0.98088\\
\hline
A.7 &  0.00037 & 0.96472 & 0.96472 & 0.96472 & 0.96472\\
\hline
A.8 &  0.00714 & 0.93729 & 0.93729 & 0.93729 & 0.93729\\
\hline
A.9 & 0.04889 & 0.99073 & 0.99073 & 0.99073 & 0.99073\\
\hline
A.10 & 0.0751  & 0.99137 & 0.99137 & 0.99137 & 0.99137\\
 \hline
Average & 0.091967 & 0.92544 & 0.9118 & 0.92766 & 0.92545\\
 \hline
Sd & 0.14025 & 0.10266 & 0.18927 & 0.15079 & 0.15035\\
 \hline
SEM & 0.044352 & 0.032465 & 0.059853 & 0.047683 & 0.047546\\
 \hline
p-value & 0.067972 & 3.9164e-10 & 9.861e-08 & 1.1585e-08 & 1.1534e-08\\
 \hline
95\% CI & (-0.00836,0.1923) & (0.852,0.99888) & (0.7764,1.0472) & (0.81979,1.0355) & 	(0.8179,1.033)\\
 \hline
\end{tabular}
}
\end{center}
\end{table}

\noindent

\begin{figure}
\centering
\includegraphics[width=0.8\textwidth]{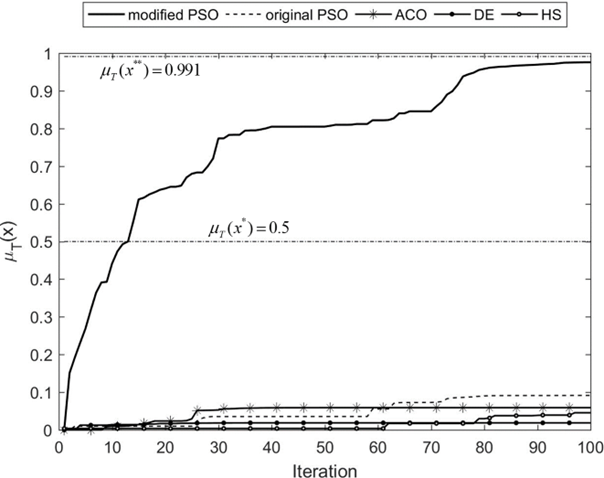}
\caption{The results (averaged over 30 runs) of the MPSO, OPSO, CACO, DE and HS algorithms on test problem A.1 and their differences from the best super-optimum found by the linearization approach.}
\end{figure}

\begin{figure}
\centering
\includegraphics[width=0.8\textwidth]{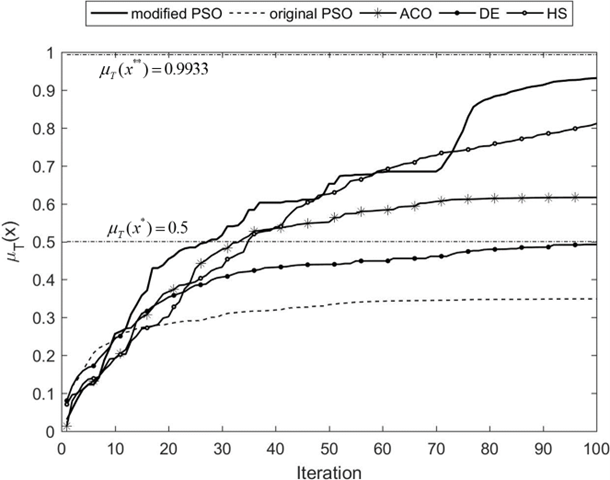}
\caption{The results (averaged over 30 runs) of the MPSO, OPSO, CACO, DE and HS algorithms on test problem A.2 and their differences from the best super-optimum found by the linearization approach.}
\end{figure}

\begin{figure}
\centering
\includegraphics[width=0.8\textwidth]{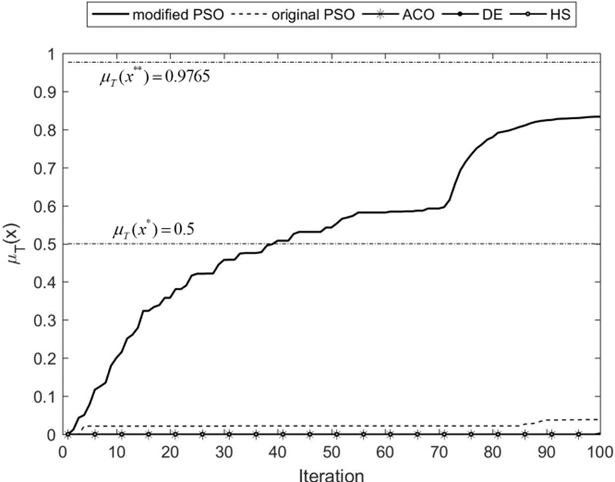}
\caption{The results (averaged over 30 runs) of the MPSO, OPSO, CACO, DE and HS algorithms on test problem A.3 and their differences from the best super-optimum found by the linearization approach.}
\end{figure}

\begin{figure}
\centering
\includegraphics[width=0.8\textwidth]{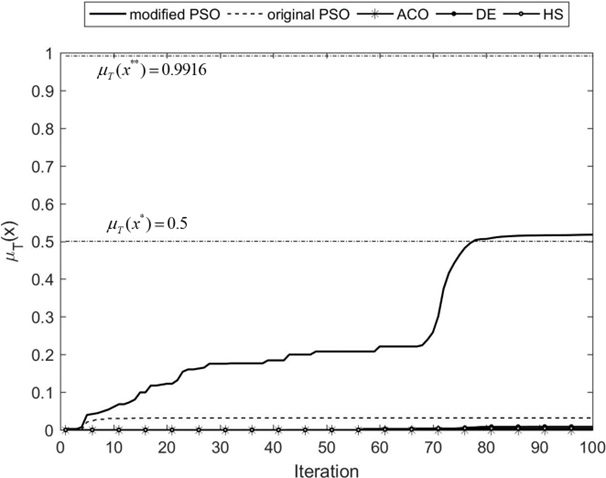}
\caption{The results (averaged over 30 runs) of the MPSO, OPSO, CACO, DE and HS algorithms on test problem A.4 and their differences from the best super-optimum found by the linearization approach.}
\end{figure}

\begin{figure}
\centering
\includegraphics[width=0.8\textwidth]{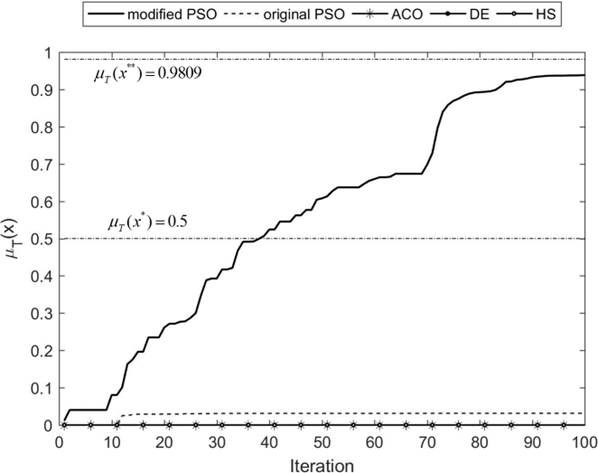}
\caption{The results (averaged over 30 runs) of the MPSO, OPSO, CACO, DE and HS algorithms on test problem A.6 and their differences from the best super-optimum found by the linearization approach.}
\end{figure}

\begin{figure}
\centering
\includegraphics[width=0.8\textwidth]{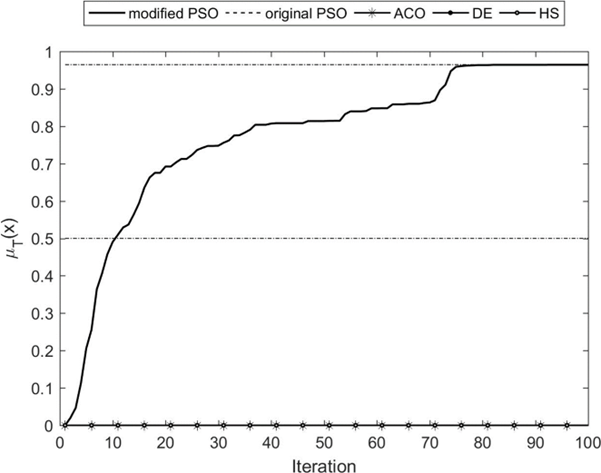}
\caption{The results (averaged over 30 runs) of the MPSO, OPSO, CACO, DE and HS algorithms on test problem A.7 and their differences from the best super-optimum found by the linearization approach.}
\end{figure}

\begin{figure}
\centering
\includegraphics[width=0.8\textwidth]{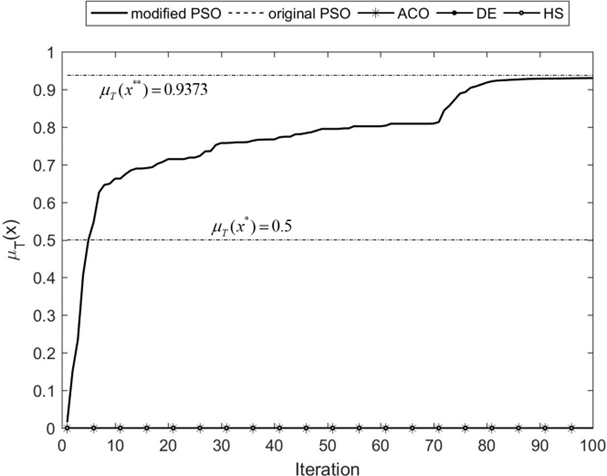}
\caption{The results (averaged over 30 runs) of the MPSO, OPSO, CACO, DE and HS algorithms on test problem A.8 and their differences from the best super-optimum found by the linearization approach.}
\end{figure}

\noindent \\ \\ \\ \\ \\ \\ \\

\vskip 0.2in
\noindent \textbf{Conclusion }
\vskip 0.2in
\noindent In this paper, a new algorithm was presented for solving the fuzzy relational inequalities with fuzzy constraints (FRI-FC) defined by the product t-norm. The linear optimization problem with respect to the regions formed as FRI-FC was studied and a set of equivalent optimization problems was introduced by the theoretical aspects of the fuzzy relational inequalities. It is shown that the main problem can be converted into an equivalent linear model that can be solved by many efficient methods such as the simplex algorithm. Moreover, three simplification operations were presented to convert the problem into a more simplified one. A comparison was made between the proposed method and the modified PSO which solves the linear optimization problems subjected to the FRI-FC regions. Furthermore, since the main problem can be transformed into an unconstrained optimization problem, the solutions of the linearization method and the modified PSO were compared to the results obtained by some well-known heuristic algorithms such as original PSO (OPSO), continuous Ant Colony Optimization (CACO), Differential Evolution (DE) and Harmony Search (HS).  On the contrary to the OPSO, CACO, DE and HS algorithms, the solutions found by the modified PSO were mostly very close to the best super-optima generated by the linearization method. From the results of the errors and times of the execution of the proposed method, we can observe that our algorithm finds the best super-optima with the least constraints violation and finishes an experiment in less than 0.996 seconds.  As future works, we aim at testing our algorithm in other types of FRI-FC problems.

\noindent \textbf{}

\vskip 0.2in
\noindent \textbf{Appendix A}
\vskip 0.2in

\noindent \textbf{Test Problem A.1}
\vskip 0.1in
\noindent \(c^{T} =[
    \begin{array}{l}
        {-1.7005 \; \; -4.3370 \; \; -3.5848 \; \; 7.7951 \; \; 0.4787 \; \; 9.3360}
    \end{array}
]\)
\vskip 0.1in
\noindent \(A=\left[
    \begin{array}{l}
        {0.1359 \; \; 0.8372 \; \; 0.1439 \; \; 0.8102 \; \; 0.8317 \; \; 0.0801} \\
        {0.0866 \; \; 0.1400 \; \; 0.9757 \; \; 0.1262 \; \; 0.7061 \; \; 0.8810} \\
        {0.7325 \; \; 0.0336 \; \; 0.3292 \; \; 0.6075 \; \; 0.6609 \; \; 0.6825} \\
        {0.8851 \; \; 0.3458 \; \; 0.6788 \; \; 0.5612 \; \; 0.8223 \; \; 0.3284}
    \end{array}
\right]\)
\vskip 0.1in
\noindent \(b^{T} =[
    \begin{array}{l}
        {0.2014 \; \; 0.0161 \; \; 0.6792 \; \; 0.8360}
    \end{array}
]\)

\vskip 0.2in
\noindent \textbf{Test Problem A.2}
\vskip 0.1in
\noindent \(c^{T} =[
    \begin{array}{l}
        {6.1322 \; \; -1.5004 \; \; 4.2041 \; \; 0.1303 \; \; -7.0502 \; \; -9.3533}
    \end{array}
]\)
\vskip 0.1in
\noindent \(A=\left[
    \begin{array}{l}
        {0.5847 \; \; 0.1338 \; \; 0.4806 \; \; 0.2675 \; \; 0.3970 \; \; 0.3215} \\
        {0.1819 \; \; 0.1038 \; \; 0.0074 \; \; 0.0476 \; \; 0.1761 \; \; 0.0927} \\
        {0.2044 \; \; 0.2334 \; \; 0.1171 \; \; 0.3456 \; \; 0.2832 \; \; 0.3481} \\
        {0.0659 \; \; 0.4362 \; \; 0.5850 \; \; 0.3090 \; \; 0.7723 \; \; 0.4813} \\
        {0.8992 \; \; 0.0619 \; \; 0.7021 \; \; 0.0030 \; \; 0.8549 \; \; 0.8145}
    \end{array}
\right]\)
\vskip 0.1in
\noindent \(b^{T} =[
    \begin{array}{l}
        {0.9366 \; \; 0.1139 \; \; 0.2031 \; \; 0.8282 \; \; 0.8752}
    \end{array}
]\)

\vskip 0.2in
\noindent \textbf{Test Problem A.3}
\vskip 0.1in
\noindent \(c^{T} =[
    \begin{array}{l}
        {6.8552 \; \; -1.8291 \; \; 8.6740 \; \; -6.4465 \; \; 5.9379 \; \; -1.7184}
    \end{array}
]\)
\vskip 0.1in
\noindent \(A=\left[
    \begin{array}{l}
        {0.3248 \; \; 0.1426 \; \; 01232 \; \; 0.8051 \; \; 0.1083 \; \; 0.2499} \\
        {0.7076 \; \; 0.2252 \; \; 0.7964 \; \; 0.0571 \; \; 0.6141 \; \; 0.0879} \\
        {0.1289 \; \; 0.7975 \; \; 0.0480 \; \; 0.0046 \; \; 0.0144 \; \; 0.1980} \\
        {0.5316 \; \; 0.8095 \; \; 0.8496 \; \; 0.1967 \; \; 0.5240 \; \; 0.0389} \\
        {0.6753 \; \; 0.6560 \; \; 0.3397 \; \; 0.4811 \; \; 0.3146 \; \; 0.4192} \\
        {0.0151 \; \; 0.7319 \; \; 0.4089 \; \; 0.6728 \; \; 0.6220 \; \; 0.7515}
    \end{array}
\right]\)
\vskip 0.1in
\noindent \(b^{T} =[
    \begin{array}{l}
        {0.6034 \; \; 0.4401 \; \; 0.5971 \; \; 0.1162 \; \; 0.0642 \; \; 0.8811}
    \end{array}
]\)

\vskip 0.2in
\noindent \textbf{Test Problem A.4}
\vskip 0.1in
\noindent \(c^{T} =[
    \begin{array}{l}
        {-5.9968 \; \; 1.6453 \; \; 7.2174 \; \; -6.3971 \; \; 2.9091 \; \; 1.7317 \; \; -9.5134}
    \end{array}
]\)
\vskip 0.1in
\noindent \(A=\left[
    \begin{array}{l}
        {0.1556 \; \; 0.7437 \; \; 0.3357 \; \; 0.9491 \; \; 0.8195 \; \; 0.7906 \; \; 0.0144} \\
        {0.8460 \; \; 0.1222 \; \; 0.4980 \; \; 0.7749 \; \; 0.7375 \; \; 0.4606 \; \; 0.3457} \\
        {0.6406 \; \; 0.5979 \; \; 0.9258 \; \; 0.2752 \; \; 0.2663 \; \; 0.4620 \; \; 0.6144} \\
        {0.5399 \; \; 0.3409 \; \; 0.4983 \; \; 0.9553 \; \; 0.7804 \; \; 0.4636 \; \; 0.7174} \\
        {0.1448 \; \; 0.1067 \; \; 0.8405 \; \; 0.8741 \; \; 0.6756 \; \; 0.8831 \; \; 0.6736}
    \end{array}
\right]\)
\vskip 0.1in
\noindent \(b^{T} =[
    \begin{array}{l}
        {0.0322 \; \; 0.4687 \; \; 0.8801 \; \; 0.6245 \; \; 0.8080}
    \end{array}
]\)

\vskip 0.2in
\noindent \textbf{Test Problem A.5}
\vskip 0.1in
\noindent \(c^{T} =[
    \begin{array}{l}
        {-2.2349 \; \; 9.8684 \; \; 6.0145 \; \; -0.1409 \; \; 4.4755 \; \; 5.1849 \; \; -5.2513 \; \; 7.4044}
    \end{array}
]\)
\vskip 0.1in
\noindent \(A=\left[
    \begin{array}{l}
        {0.6685 \; \; 0.8277 \; \; 0.7281 \; \; 0.1500 \; \; 0.0963 \; \; 0.9766 \; \; 0.0665 \; \; 0.0789} \\
        {0.1002 \; \; 0.7566 \; \; 0.3847 \; \; 0.2817 \; \; 0.6514 \; \; 0.7538 \; \; 0.2835 \; \; 0.4277} \\
        {0.5461 \; \; 0.8736 \; \; 0.7027 \; \; 0.3276 \; \; 0.9032 \; \; 0.3744 \; \; 0.5540 \; \; 0.5274} \\
        {0.5924 \; \; 0.2766 \; \; 0.5094 \; \; 0.3399 \; \; 0.3381 \; \; 0.9006 \; \; 0.5930 \; \; 0.4494} \\
        {0.9426 \; \; 0.0540 \; \; 0.7771 \; \; 0.8831 \; \; 0.1892 \; \; 0.4538 \; \; 0.4534 \; \; 0.3912} \\
        {0.3466 \; \; 0.2082 \; \; 0.2217 \; \; 0.2299 \; \; 0.7907 \; \; 0.0611 \; \; 0.3830 \; \; 0.9836}
    \end{array}
\right]\)
\vskip 0.1in
\noindent \(b^{T} =[
    \begin{array}{l}
        {0.9642 \; \; 0.0346 \; \; 0.1792 \; \; 0.7077 \; \; 0.9931 \; \; 0.3542}
    \end{array}
]\)

\vskip 0.2in
\noindent \textbf{Test Problem A.6}
\vskip 0.1in
\noindent \(c^{T} =[
    \begin{array}{l}
        {3.9815 \; \; 7.8181 \; \; 9.1858 \; \; 0.9443 \; \; -7.2275 \; \; -7.0141 \; \; -4.8498}
    \end{array}
]\)
\vskip 0.1in
\noindent \(A=\left[
    \begin{array}{l}
        {0.8407 \; \; 0.2511 \; \; 0.9172 \; \; 0.0540 \; \; 0.0119 \; \; 0.6020 \; \; 0.2290} \\
        {0.2543 \; \; 0.6160 \; \; 0.2858 \; \; 0.5308 \; \; 0.3371 \; \; 0.2630 \; \; 0.9133} \\
        {0.8143 \; \; 0.4733 \; \; 0.7572 \; \; 0.7792 \; \; 0.1622 \; \; 0..6541 \; \; 0.1524} \\
        {0.2435 \; \; 0.3517 \; \; 0.7537 \; \; 0.9340 \; \; 0.7943 \; \; 0.6892 \; \; 0.8258} \\
        {0.9293 \; \; 0.8308 \; \; 0.3804 \; \; 0.1299 \; \; 0.3112 \; \; 0.7482 \; \; 0.5383} \\
        {0.3500 \; \; 0.5853 \; \; 0.5678 \; \; 0.5688 \; \; 0.5285 \; \; 0.4505 \; \; 0.9961} \\
        {0.1966 \; \; 0.5497 \; \; 0.0759 \; \; 0.4694 \; \; 0.1656 \; \; 0.0838 \; \; 0.0782}
    \end{array}
\right]\)
\vskip 0.1in
\noindent \(b^{T} =[
    \begin{array}{l}
        {0.4427 \; \; 0.1067 \; \; 0.9619 \; \; 0.0046 \; \; 0.7749 \; \; 0.8173 \; \; 0.8687}
    \end{array}
]\)

\vskip 0.2in
\noindent \textbf{Test Problem A.7}
\vskip 0.1in
\noindent \(c^{T} =[
    \begin{array}{l}
        {-6.2250 \; \; 7.3391 \; \; 7.0215 \; \; 6.5683 \; \; 1.2437 \; \; 5.2404 \; \; 4.1114 \; \; 1.4945 \; \; 7.7537 \; \; -0.0706}
    \end{array}
]\)
\vskip 0.1in
\noindent \(A=\left[
    \begin{array}{l}
        {0.8376 \; \; 0.5607 \; \; 0.4054 \; \; 0.9161 \; \; 0.7397 \; \; 0.0716 \; \; 0.3919 \; \; 0.2039 \; \; 0.2089 \; \; 0.5333} \\
        {0.0516 \; \; 0.9841 \; \; 0.0849 \; \; 0.2607 \; \; 0.1426 \; \; 0.3386 \; \; 0.1199 \; \; 0.4257 \; \; 0.5544 \; \; 0.3795} \\
        {0.4259 \; \; 0.7705 \; \; 0.8464 \; \; 0.3001 \; \; 0.4591 \; \; 0.6629 \; \; 0.3465 \; \; 0.5807 \; \; 0.0754 \; \; 0.7783} \\
        {0.2721 \; \; 0.0729 \; \; 0.5261 \; \; 0.4080 \; \; 0.6138 \; \; 0.5058 \; \; 0.5298 \; \; 0.0238 \; \; 0.5426 \; \; 0.8281} \\
        {0.7275 \; \; 0.3988 \; \; 0.5304 \; \; 0.8283 \; \; 0.3916 \; \; 0.9934 \; \; 0.1053 \; \; 0.7491 \; \; 0.2099 \; \; 0.8877} \\
        {0.8102 \; \; 0.0718 \; \; 0.8724 \; \; 0.4534 \; \; 0.4194 \; \; 0.4759 \; \; 0.3741 \; \; 0.0360 \; \; 0.7298 \; \; 0.2837} \\
        {0.4767 \; \; 0.1241 \; \; 0.7166 \; \; 0.9912 \; \; 0.3778 \; \; 0.8198 \; \; 0.6201 \; \; 0.2147 \; \; 0.1651 \; \; 0.6314}
    \end{array}
\right]\)
\vskip 0.1in
\noindent \(b^{T} =[
    \begin{array}{l}
        {0.1699 \; \; 0.3880 \; \; 0.1850 \; \; 0.2983 \; \; 0.0652 \; \; 0.7516 \; \; 0.0027}
    \end{array}
]\)

\vskip 0.2in
\noindent \textbf{Test Problem A.8}
\vskip 0.1in
\noindent \(c^{T} =[
    \begin{array}{l}
        {-0.6133 \; \; 9.5958 \; \; 8.9277 \; \; 2.5328 \; \; 4.5288 \; \; -3.2146 \; \; 8.4352 \; \; 9.0807 \; \; 9.8011 \; \; 2.4952}
    \end{array}
]\)
\vskip 0.1in
\noindent \(A=\left[
    \begin{array}{l}
        {0.9576 \; \; 0.1524 \; \; 0.4591 \; \; 0.1542 \; \; 0.6114 \; \; 0.3348 \; \; 0.5124 \; \; 0.5708 \; \; 0.5980 \; \; 0.4008} \\
        {0.0317 \; \; 0.1542 \; \; 0.0775 \; \; 0.4483 \; \; 0.4208 \; \; 0.2367 \; \; 0.2025 \; \; 0.2892 \; \; 0.7608 \; \; 0.7968} \\
        {0.4176 \; \; 0.9340 \; \; 0.5365 \; \; 0.5206 \; \; 0.3109 \; \; 0.6503 \; \; 0.9856 \; \; 0.7474 \; \; 0.0019 \; \; 0.3848} \\
        {0.5935 \; \; 0.2486 \; \; 0.9793 \; \; 0.3124 \; \; 0.2278 \; \; 0.9957 \; \; 0.1159 \; \; 0.8994 \; \; 0.6739 \; \; 0.7614} \\
        {0.2106 \; \; 0.3407 \; \; 0.2622 \; \; 0.4472 \; \; 0.3391 \; \; 0.7915 \; \; 0.3352 \; \; 0.6066 \; \; 0.1036 \; \; 0.2750} \\
        {0.2554 \; \; 0.2953 \; \; 0.2070 \; \; 0.3073 \; \; 0.8167 \; \; 0.2889 \; \; 0.1124 \; \; 0.7823 \; \; 0.5541 \; \; 0.2769} \\
        {0.6388 \; \; 0.2688 \; \; 0.6292 \; \; 0.7254 \; \; 0.9561 \; \; 0.5072 \; \; 0.0170 \; \; 0.4677 \; \; 0.3709 \; \; 0.7724} \\
        {0.5010 \; \; 0.2386 \; \; 0.5260 \; \; 0.3181 \; \; 0.2825 \; \; 0.2465 \; \; 0.6896 \; \; 0.1970 \; \; 0.4726 \; \; 0.8429}
    \end{array}
\right]\)
\vskip 0.1in
\noindent \(b^{T} =[
    \begin{array}{l}
        {0.9750 \; \; 0.0237 \; \; 0.5693 \; \; 0.5947 \; \; 0.0129 \; \; 0.7476 \; \; 0.8917 \; \; 0.9195}
    \end{array}
]\)

\vskip 0.2in
\noindent \textbf{Test Problem A.9}
\vskip 0.1in
\noindent \(c^{T} =[
    \begin{array}{l}
        {2.7276 \; \; -9.1572 \; \; -7.1452 \; \; 0.9314 \; \; -1.9806 \; \; -0.2452 \; \; 8.1128 \; \; 2.4546 \; \; 2.85}
    \end{array}
]\)
\vskip 0.1in
\noindent \(A=\left[
    \begin{array}{l}
        {0.6637 \; \; 0.7432 \; \; 0.8188 \; \; 0.5660 \; \; 0.6163 \; \; 0.6756 \; \; 0.8136 \; \; 0.9008 \; \; 0.1373} \\
        {0.4567 \; \; 0.8137 \; \; 0.1834 \; \; 0.5090 \; \; 0.8100 \; \; 0.8336 \; \; 0.0360 \; \; 0.6520 \; \; 0.7940} \\
        {0.9142 \; \; 0.6858 \; \; 0.3362 \; \; 0.4732 \; \; 0.9350 \; \; 0.3233 \; \; 0.1749 \; \; 0.3490 \; \; 0.6930} \\
        {0.3936 \; \; 0.8700 \; \; 0.9014 \; \; 0.4516 \; \; 0.8772 \; \; 0.2706 \; \; 0.2962 \; \; 0.7862 \; \; 0.8683} \\
        {0..5089 \; \; 0.8362 \; \; 0.9459 \; \; 0.6415 \; \; 0.7242 \; \; 0.3912 \; \; 0.0342 \; \; 0.8280 \; \; 0.1411} \\
        {0.1390 \; \; 0.9381 \; \; 0.4557 \; \; 0.4086 \; \; 0.4804 \; \; 0.7904 \; \; 0.2785 \; \; 0.5207 \; \; 0.4494} \\
        {0.1239 \; \; 0.4953 \; \; 0.3745 \; \; 0.3644 \; \; 0.8924 \; \; 0.5274 \; \; 0.5470 \; \; 0.6684 \; \; 0.5228} \\
        {0.6215 \; \; 0.6577 \; \; 0.2243 \; \; 0.2977 \; \; 0.8589 \; \; 0.8018 \; \; 0.0108 \; \; 0.4651 \; \; 0.9415} \\
        {0.6028 \; \; 0.2612 \; \; 0.0686 \; \; 0.7647 \; \; 0.0017 \; \; 0.2003 \; \; 0.2626 \; \; 0.1225 \; \; 05671}
    \end{array}
\right]\)
\vskip 0.1in
\noindent \(b^{T} =[
    \begin{array}{l}
        {0.6412 \; \; 0.0190 \; \; 0.1702 \; \; 0.6046 \; \; 0.2470 \; \; 0.8351 \; \; 0.7981 \; \; 0.4645 \; \; 0.6098}
    \end{array}
]\)

\vskip 0.2in
\noindent \textbf{Test Problem A.10}
\vskip 0.1in
\noindent \(c^{T} =[
    \begin{array}{l}
        {0.1117 \; \; -5.4486 \; \; 4.4966 \; \; 2.1297 \; \; -0.4456 \; \; -6.9081 \; \; 6.4627 \; \; 4.4522 \; \; 5.6057 \; \; -5.7421}
    \end{array}
]\)
\vskip 0.1in
\noindent \(A=\left[
    \begin{array}{l}
        {0.7906 \; \; 0.9957 \; \; 0.0597 \; \; 0.5569 \; \; 0.8766 \; \; 0.4642 \; \; 0.7950 \; \; 0.4590 \; \; 0.5371 \; \; 0.1593} \\
        {0.4117 \; \; 0.3598 \; \; 0.0890 \; \; 0.7960 \; \; 0.6433 \; \; 0.7929 \; \; 0.9148 \; \; 0.2328 \; \; 0.5097 \; \; 0.8664} \\
        {0.8080 \; \; 0.8941 \; \; 0.6725 \; \; 0.9067 \; \; 0.6509 \; \; 0.5584 \; \; 0.0008 \; \; 0.0250 \; \; 0.5247 \; \; 0.0786} \\
        {0.0878 \; \; 0.9908 \; \; 0.1796 \; \; 0.9319 \; \; 0.1598 \; \; 0.0141 \; \; 0.8120 \; \; 0.0548 \; \; 0.7654 \; \; 0.8247} \\
        {0.0338 \; \; 0.7684 \; \; 0.8985 \; \; 0.1536 \; \; 0.8718 \; \; 0.2336 \; \; 0.8460 \; \; 0.5478 \; \; 0.6549 \; \; 0.7402} \\
        {0.8519 \; \; 0.1464 \; \; 0.0081 \; \; 0.5379 \; \; 0.8072 \; \; 0.3497 \; \; 0.9215 \; \; 0.9925 \; \; 0.1810 \; \; 0.1014} \\
        {0.2547 \; \; 0.9694 \; \; 0.0436 \; \; 0.7910 \; \; 0.0878 \; \; 0.7878 \; \; 0.4435 \; \; 0.6274 \; \; 0.8116 \; \; 0.3670} \\
        {0.1954 \; \; 0.8315 \; \; 0.7230 \; \; 0.8740 \; \; 0.4403 \; \; 0.2030 \;  0.0607 \; \; 0.2924 \; \; 0.2905 \; \; 0.8977} \\
        {0.7841 \; \; 0.2572 \; \; 0.5130 \; \; 0.9844 \; \; 0.5482 \; \; 0.5547 \; \; 0.2540 \; \; 0.7913 \; \; 0.8254 \; \; 0.0586} \\
        {0.1628 \; \; 0.6455 \; \; 0.8100 \; \; 0.2105 \; \; 0.3994 \; \; 0.9538 \; \; 0.9343 \; \; 0.3079 \; \; 0.1052 \; \; 0.5025}
    \end{array}
\right]\)
\vskip 0.1in
\noindent \(b^{T} =[
    \begin{array}{l}
        {0.0294 \; \; 0.9041  \; \; 0.9020 \; \; 0.9824 \; \; 0.8485 \; \; 0.1962 \; \; 0.4210 \; \; 0.2162 \; \; 0.8787 \; \; 0.3932}
    \end{array}
]\)

\vskip 0.2in
\noindent \textbf{Acknowledgment}
\vskip 0.2in
\noindent We are very grateful to the anonymous referees for their comments and suggestions, which were very helpful in improving the paper.

\vskip 0.2in
\noindent \textbf{References}
\vskip 0.2in
\noindent [1]. C. W. Chang, B. S. Shieh, Linear optimization problem constrained by fuzzy max–min relation equations, Information Sciences 234 (2013) 71–79.

\noindent [2]. L. Chen, P. P. Wang, Fuzzy relation equations (i): the general and specialized solving algorithms, Soft Computing 6 (5) (2002) 428-435.

\noindent [3]. L. Chen, P. P. Wang, Fuzzy relation equations (ii): the branch-point-solutions and the categorized minimal solutions, Soft Computing 11 (1) (2007) 33-40.

\noindent [4]. S. Dempe, A. Ruziyeva, On the calculation of a membership function for the solution of a fuzzy linear optimization problem, Fuzzy Sets and Systems 188 (2012) 58-67.

\noindent [5]. A. Di Nola, S. Sessa, W. Pedrycz, E. Sanchez, Fuzzy relational Equations and their applications in knowledge engineering, Dordrecht: Kluwer Academic Press, 1989.

\noindent [6]. D. Dubey, S. Chandra, A. Mehra, Fuzzy linear programming under interval uncertainty based on IFS representation, Fuzzy Sets and Systems 188 (2012) 68-87.

\noindent [7]. D. Dubois, H. Prade, Fundamentals of Fuzzy Sets, Kluwer, Boston, 2000.

\noindent [8]. Y. R. Fan, G. H. Huang, A. L. Yang, Generalized fuzzy linear programming for decision making under uncertainty: Feasibility of fuzzy solutions and solving approach, Information Sciences 241 (2013) 12-27.

\noindent [9]. S.C. Fang, G. Li, Solving fuzzy relational equations with a linear objective function, Fuzzy Sets and Systems 103 (1999) 107-113.

\noindent [10]. S. Freson, B. De Baets, H. De Meyer, Linear optimization with bipolar max–min constraints, Information Sciences 234 (2013) 3–15.

\noindent [11]. Z.W. Geem, J.H. Kim, G.V. Loganathan, A new heuristic optimization algorithm: harmony search, Simulation 76(2)(2001) 60-68. 

\noindent [12]. A. Ghodousian, E. Khorram, Fuzzy linear optimization in the presence of the fuzzy relation inequality constraints with max-min composition, Information Sciences 178 (2008) 501-519.

\noindent [13]. A. Ghodousian, E. Khorram, Linear optimization with an arbitrary fuzzy relational inequality, Fuzzy Sets and Systems 206 (2012) 89-102.

\noindent [14]. A. Ghodousian, M. Raeisian Parvari, A modiﬁed PSO algorithm for linear optimization problem subject to the generalized fuzzy relational inequalities with fuzzy constraints (FRI-FC), Information Sciences 418–419 (2017) 317–345.

\noindent [15].  A. Ghodousian, A. Babalhavaeji, An efﬁcient genetic algorithm for solving nonlinear optimization problems deﬁned with fuzzy relational equations and max-Lukasiewicz composition, Applied Soft Computing 69 (2018) 475–492.

\noindent [16]. A. Ghodousian, M. Naeeimib, A. Babalhavaeji, Nonlinear optimization problem subjected to fuzzy relational equations deﬁned by Dubois-Prade family of t-norms, Computers \& Industrial Engineering 119 (2018) 167–180.

\noindent [17]. A. Ghodousian, Optimization of linear problems subjected to the intersection of two fuzzy relational inequalities defined by Dubois-Prade family of t-norms, Information Sciences 503 (2019) 291–306.

\noindent [18]. A. Ghodousian, An algorithm for solving linear optimization problems subjected to the intersection of two fuzzy relational inequalities defined by Frank family of t-norms, International Journal in Foundations of Computer Science and Technology 8(3)(2018) 1-20.

\noindent [19]. A. Ghodousian, M. Raeisian Parvari, R. Rabie, T. Azarnejad, Solving a Non-Linear Optimization Problem Constrained by a Non-Convex Region Defined by Fuzzy Relational Equations and Schweizer-Sklar Family of T-Norms, American Journal of Computation, Communication and Control 5(2)(2018) 68-87. 

\noindent [20]. F. Guo, Z. Q. Xia, An algorithm for solving optimization problems with one linear objective function and finitely many constraints of fuzzy relation inequalities, Fuzzy Optimization and Decision Making 5 (2006) 33-47.

\noindent [21]. F. F. Guo, L. P. Pang, D. Meng, Z. Q. Xia, An algorithm for solving optimization problems with fuzzy relational inequality constraints, Information Sciences 252 ( 2013) 20-31.

\noindent [22]. S. M. Guu, Y. K. Wu, Minimizing a linear objective function under a max-t-norm fuzzy relational equation constraint, Fuzzy Sets and Systems 161 (2010) 285-297.

\noindent [23]. S. M. Guu, Y. K. Wu, Minimizing a linear objective function with fuzzy relation equation constraints, Fuzzy Optimization and Decision Making 12 (2002) 1568-4539.

\noindent [24]. J. Kennedy, R. C. Eberhart, Particle Swarm Optimization, in: Proceedings of IEEE International Conference on Neural Networks, vol. 4, 1995, pp.1942-1948.

\noindent [25]. H. C. Lee, S. M. Guu, On the optimal three-tier multimedia streaming services, Fuzzy Optimization and Decision Making 2(1) (2002) 31-39.

\noindent [26]. P. K. Li, S. C. Fang, On the resolution and optimization of a system of fuzzy relational equations with sup-t composition, Fuzzy Optimization and Decision Making 7 (2008) 169-214.

\noindent [27]. J. X. Li, S. J. Yang, Fuzzy relation inequalities about the data transmission mechanism in bittorrent-like peer-to-peer file sharing systems, in: Proceedings of the 9th International Conference on Fuzzy Systems and Knowledge discovery (FSKD 2012), pp. 452-456.  

\noindent [28]. P. Li, Y. Liu, Linear optimization with bipolar fuzzy relational equation constraints using lukasiewicz triangular norm, Soft Computing 18 (2014) 1399-1404.

\noindent [29]. J. L. Lin, Y. K. Wu, S. M. Guu, On fuzzy relational equations and the covering problem, Information Sciences 181 (2011) 2951-2963.

\noindent [30]. J. L. Lin, On the relation between fuzzy max-archimedean t-norm relational equations and the covering problem, Fuzzy Sets and Systems 160 (2009) 2328-2344. 

\noindent [31]. C. C. Liu, Y. Y. Lur, Y. K. Wu, Linear optimization of bipolar fuzzy relational equations with max-Łukasiewicz composition, Information Sciences 360 (2016) 149–162.

\noindent [32]. J. Loetamonphong, S. C. Fang, Optimization of fuzzy relation equations with max-product composition, Fuzzy Sets and Systems 118 (2001) 509-517.

\noindent [33]. A. V. Markovskii, On the relation between equations with max-product composition and the covering problem, Fuzzy Sets and Systems 153 (2005) 261-273.

\noindent [34]. M. Mizumoto, H. J. Zimmermann, Comparison of fuzzy reasoning method, Fuzzy Sets and Systems 8 (1982) 253-283.

\noindent [35]. W. Pedrycz, Granular Computing: Analysis and Design of Intelligent Systems, CRC Press, Boca Raton, 2013.

\noindent [36]. W. Pedrycz, Fuzzy relational equations with generalized connectives and their applications, Fuzzy Sets and Systems 10 (1983) 185-201.

\noindent [37]. W. Pedrycz, Solving fuzzy relational equations through logical filtering, Fuzzy Sets and Systems 81 (1996) 355-363.

\noindent [38]. K.V. Price, R.M. Storn, J.A. Lampinen, Differential evolution-A practical approach to global optimization, Springer, Verlag Berlin Heidelberg, 2005.

\noindent [39]. X. B. Qu, X. P. Wang, Man-hua. H. Lei, Conditions under which the solution sets of fuzzy relational equations over complete Brouwerian lattices form lattices, Fuzzy Sets and Systems 234 (2014) 34-45.

\noindent [40]. X. B. Qu, X. P. Wang, Minimization of linear objective functions under the constraints expressed by a system of fuzzy relation equations, Information Sciences 178 (2008) 3482-3490.

\noindent [41]. E. Sanchez, Solution in composite fuzzy relation equations: application to medical diagnosis in Brouwerian logic, in: M.M. Gupta. G.N. Saridis, B.R. Games (Eds.), Fuzzy Automata and Decision Processes, North-Holland, New York, 1977, pp. 221-234.

\noindent [42]. E. Sanchez, Resolution of eigen fuzzy sets equations, Fuzzy Sets and Systems 1 (1978) 69-75.

\noindent [43]. B. S. Shieh, Infinite fuzzy relation equations with continuous t-norms, Information Sciences 178 (2008) 1961-1967.

\noindent [44]. B. S. Shieh, Minimizing a linear objective function under a fuzzy max-t-norm relation equation constraint, Information Sciences 181 (2011) 832-841.

\noindent [45]. K. Socha, M. Dorigo, “Ant colony optimization for continuous domain”, European Journal of operational Research 185 (2008) 1155-1173.

\noindent [46]. F. Sun, Conditions for the existence of the least solution and minimal solutions to fuzzy relation equations over complete Brouwerian lattices, Information Sciences 205 (2012) 86-92.

\noindent [47]. F. Sun, X. P. Wang, x. B. Qu, Minimal join decompositions and their applications to fuzzy relation equations over complete Brouwerian lattices, Information Sciences 224 (2013) 143-151. 

\noindent [48]. Y. K. Wu, Optimization of fuzzy relational equations with max-av composition, Information Sciences 177 (2007) 4216-4229.

\noindent [49]. Y. K. Wu, S. M. Guu, Minimizing a linear function under a fuzzy max-min relational equation constraints, Fuzzy Sets and Systems 150 (2005) 147-162.

\noindent [50]. Y. K. Wu, S. M. Guu, An efficient procedure for solving a fuzzy relation equation with max-Archimedean t-norm composition, IEEE Transactions on Fuzzy Systems 16 (2008) 73-84.

\noindent [51]. Y. K. Wu, S. M. Guu, J. Y. Liu, Reducing the search space of a linear fractional programming problem under fuzzy relational equations with max-Archimedean t-norm composition, Fuzzy Sets and Systems 159 (2008) 3347-3359.

\noindent [52]. Q. Q. Xiong, X. P. Wang, Fuzzy relational equations on complete Brouwerian lattices, Information Sciences 193 (2012) 141-152.

\noindent [53]. X. P. Yang, X. G. Zhou, B. Y. Cao, Latticized linear programming subject to max-product fuzzy relation inequalities with application in wireless communication, Information Sciences 358–359 (2016) 44–55.

\noindent [54]. S. J. Yang, An algorithm for minimizing a linear objective function subject to the fuzzy relation inequalities with addition-min composition, Fuzzy Sets and Systems 255 (2014) 41-51.  \textbf{}

\end{document}